\documentclass[11pt, 
]{amsart}
\usepackage{ 
amssymb 
} 
\newcommand{\no}[1]{ }

\date{\today}

\newtheorem{theorem}{Theorem}[section] 
\newtheorem{proposition}{Proposition}[section] 
\newtheorem{lemma}{Lemma}[section] 
\newtheorem{definition}{Definition}[section] 
\newtheorem{corollary}{Corollary}[section]

\theoremstyle{remark} 
\newtheorem{remark}{Remark}[section]

\DeclareMathOperator{\supp}{supp}

\DeclareMathOperator{\dist}{dist} 
 
\newcommand{\eps}{\varepsilon} 
\newcommand{\ep}{\epsilon}
 
\newcommand{\R}{\mathbb{R}}
\newcommand{\Id}{\mbox{Id}} 
\renewcommand{\r}[1]{(\ref{#1})} 
 
\newcommand{\be}[1]{\begin{equation}\label{#1}} 
\newcommand{\ee}{\end{equation}} 

\newcommand\scl{\mathrm{sc}}
\newcommand\Psisc{\Psi_\scl}
\newcommand\Hsc{H_\scl}
\newcommand\CI{{C}^\infty}
\newcommand\dCI{\dot{C}^\infty}
\newcommand\RR{\mathbb{R}}
\newcommand\sphere{\mathbb{S}}
\newcommand\Cx{\mathbb{C}}
\newcommand\pa{\partial}
\newcommand\diag{\mathrm{diag}}
\newcommand\foliation{\mathsf{x}}

\newcommand\tufoliation{\widetilde{\underline{\mathsf{x}}}}
\newcommand\loccoord{\mathsf{y}}
\newcommand\Foliation{\mathsf{X}}
\newcommand\Loccoord{\mathsf{Y}}
\newcommand\level{\mathsf{c}}
\newcommand\Omegaext{\hat\Omega}
\newcommand\Omegaextb{\overline{\hat\Omega}}

\newcommand\cM{\mathcal{M}}

\renewcommand{\d}{\mathrm{d}}

\newcommand{\bo}{\partial M} 

 \newcommand{\zero}{^{(0)}}
\newcommand{\mat}[4]{\left(\begin{array}{cc} #1 &#2\\#3 & #4 
\end{array}\right)}

\title[Boundary rigidity with partial data]{Boundary rigidity with partial data} 

\author[P. Stefanov]{Plamen Stefanov}
\address{Department of Mathematics, Purdue University, West Lafayette, IN 47907}
\thanks{First author partly supported by  NSF  Grant DMS-1301646}

\author[G. Uhlmann]{Gunther Uhlmann}
\address{Department of Mathematics, University of Washington, Seattle, WA 98195 and Department of Mathematics University of Helsinki, Finland}
\thanks{Second author partly supported by NSF Grant CMG-1025259 and
  DMS-1265958,  The Fondation Math\'ematiques de Paris, and a Simons felloowship}

\author[A. Vasy]{Andras Vasy}
\address{Department of Mathematics, Stanford University, Stanford  CA 94305}
\thanks{Third author partly supported by NSF Grant CMG-1025259 and
DMS-1068742.}

\usepackage{graphicx}
\graphicspath{%
    {converted_graphics/}
    {/}
}
 
\begin{document} 

\begin{abstract}
We study the boundary rigidity problem with partial data consisting  of determining locally the Riemannian metric of a Riemannian manifold with boundary from the distance function measured at pairs of points near a fixed point on the boundary. We show that one can recover uniquely and in a stable way a conformal factor near a strictly convex point where we have the information. In particular, this implies that we can determine locally the isotropic sound speed of a medium by measuring the travel times of waves joining points close to a convex point on the boundary. 

The local results lead  to a global lens rigidity uniqueness and stability result assuming that the manifold is foliated by strictly convex hypersurfaces.
\end{abstract}

\maketitle 
\section{Introduction and main results} 

Travel time tomography deals with the problem of determining the sound speed or index of refraction of a medium by measuring the travel times of waves going through the medium. This type of inverse problem, also called the inverse kinematic problem,  
arose in geophysics in an attempt to determine the
substructure of the Earth by measuring the travel times of seismic waves  at the surface.  We  consider an anisotropic index of 
refraction, i.e.,  a sound speed depending on the direction. The Earth is generally isotropic. However, more recently it has been realized, by measuring these travel times, that the inner core of the Earth exhibits anisotropic behavior  with the fast direction parallel to the Earth's spin axis, see \cite{Cre}.  In the human body, muscle tissue is anisotropic. 
As a    model of anisotropy, we consider a Riemannian metric $g=(g_{ij})$. The problem is to determine the metric from the lengths of geodesics joining points on the
boundary.

This leads to the general question of whether given a compact Riemannian
manifold with boundary $(M,g)$ one can determine the Riemannian metric in the
interior knowing the boundary distance function joining points on the boundary
$d_g(x,y)$, with $x,y\in \partial M$. This is known as the {\sl boundary rigidity problem}.  Of course, isometries preserve distance, so that the boundary rigidity problem is whether two
metrics that have the same boundary distance function are the same up to
isometry fixing the boundary.  Examples can be given of manifolds that are not boundary rigid.
Such examples show that the boundary rigidity problem should be considered under some restrictions on the geometry of the manifold. The most usual of such restrictions is simplicity of the metric.  A Riemannian manifold $(M,g)$ (or the metric $g$) is called {\sl simple} if the boundary $\partial M$ is strictly convex  (positive second fundamental form)  
and any two points $x,y \in M$ are joined by a unique minimizing geodesic. Michel conjectured \cite{Michel}
that every simple compact Riemannian manifold with boundary is boundary rigid.

Simple surfaces with boundary are boundary rigid \cite{PestovU}. In higher dimensions, simple Riemannian manifolds
with boundary are boundary rigid under some a-priori constant curvature on the manifold or special symmetries
\cite{BCG}, \cite{Gr}. Several local results near the Euclidean metric are known \cite{SU-MRL}, \cite{LassasSU}, \cite{BI}. The most general result in this direction is the generic local (with respect to the metric) one proven in \cite{SU-JAMS}. 
Surveys of some of the results can be found in \cite{I}, \cite{SU-Kawai}, \cite{Croke04}.

In this paper, we consider the boundary rigidity problem in the class of  metrics conformal to a given one and with partial (local) data, that is, we know the boundary 
distance function for points on the boundary near a given point.  Partial data problems arise naturally in applications since in many cases one doesn't have access to the whole boundary. We prove the first result on the determination
of the conformal factor locally near the boundary from partial data without assuming analyticity.  We develop a novel method to attack partial data non-linear  problems that will likely have other applications.

We now describe the known results with full data on the boundary.
Let us fix the metric $g_0$ and let $c$ be a positive smooth function on the compact manifold  with boundary  $M$. The problem is whether we can determine $c$ from $d_{c^{-2} g_0}(x,y), x,y\in \partial M.$ Notice that in this case the problem is not invariant under changes of variables that are the identity at the boundary so that we expect to be able to recover $c$  under appropriate a-priori conditions. This was proven by Mukhometov in two dimensions \cite{Mu2}, and in \cite{MuRo} in higher dimensions for the case of simple metrics. Of particular importance in applications is the case of an isotropic sound speed that is when we are in a bounded domain of Euclidean space and
$g_0$ is the Euclidean metric. This is the isotropic case. This problem was considered by Herglotz \cite{Herglotz} and Wieckert and Zoeppritz \cite{WZ} for the case of a spherical symmetric sound speed. They found a formula to recover the sound speed from the boundary distance function assuming  $\frac{\d}{\d r}(\frac{r}{c(r)})>0$. Notice that this condition is equivalent  to the existence of a strictly convex foliation and is more general than simplicity, see Section~\ref{section_6}.

From now on we will call  $d$ the function $d_{{c^{-2}}g_0}$. Below, $\tilde d$ is related to $d_{{\tilde c^{-2}}g_0}$.

The partial data problem, that we will also call the {\sl local boundary rigidity problem}\footnote{It is local in the sense that $d(x,y)$ is known locally and depends on $c$ locally; the term \textit{local} has been used before to indicate that  the metric is a priori close to a fixed one.}, in this case is whether knowledge of the distance function on part of the boundary determines the sound speed $c$
locally. Given another smooth $\tilde c$, here and below we define $\tilde L$, $\tilde \ell$ and $\tilde d$ 
in the same way but related to $\tilde c$. 
We prove
the following uniqueness result:

\begin{theorem}\label{thm_1}Let $n=\dim M\ge3$,  let $c>0$, $\tilde c>0$ be smooth  and let  $\bo$ be strictly convex with respect to both $g=c^{-2}g_0$ and $\tilde g = \tilde c^{-2}g_0$ near a fixed $p\in\bo$. Let $d(p_1,p_2)= \tilde d(p_1,p_2)$ for $p_1$, $p_2$ on $\bo$ near $p$. Then $c=\tilde c$ in $M$ near $p$.
\end{theorem}

As mentioned earlier, this is the only known result for the boundary rigidity problem with partial data except in the case that the metrics are assumed to be real-analytic \cite{LassasSU}. The latter follows from determination of the jet of the metric at a convex point from the distance function known near $p.$

The boundary rigidity problem is closely connected to the lens rigidity one. To define the latter, we first introduce the manifolds $\partial_\pm SM$, defined as the sets of all vectors $(x,v)$ with $x\in\bo$, $v$ unit in  the metric $g$,  and  pointing outside/inside $M$. 
We define the \textit{scattering relation}
\be{L}
L: \partial_- SM \longrightarrow \partial_+SM 
\ee
in the following way: for each $(x,v)\in \partial_- SM$, $L(x,v)=(y,w)$, where $(y,w)$ are the exit point and  direction, if exist, of the maximal unit speed geodesic $\gamma_{x,v}$ in the metric $g$, issued from $(x,v)$. Let 
\[
\ell: \partial_-SM \longrightarrow \R\cup\infty
\]
be its length, possibly infinite. If $\ell<\infty$, we call $M$ non-trapping. The maps $(L,\ell)$ together are called \textit{lens relation} (or lens data).

The \textit{lens rigidity} problem is whether the scattering relation $L$ (and possibly, $\ell$) determine $g$ (and the topology of $M$) up to an isometry as above. The  lens rigidity problem with partial data for a sound speed is whether we can determine the speed near some $p$ from $L$ known near the unit sphere  $S_p\bo$ considered as a subset of $\partial_-SM$, i.e., for vectors with base points close to $p$ and directions pointing into $M$ close to ones tangent to $\bo$. For general metrics, we want to recover isometric copies of the metrics locally, as above. 

We assume that $\bo$ is strictly convex at $p\in\bo$ w.r.t.\ $g$. Then the  boundary rigidity and the  lens rigidity problems with partial data are equivalent: knowing $d$ near $(p,p)$ is equivalent to knowing $L$ in some neighborhood of  $S_p\partial M$. The size of that neighborhood however depends on a priori bounds of the derivatives of the metrics with which we work.  This equivalence was first noted by Michel \cite{Michel}, since the tangential gradients of $d(x,y)$ on $\bo\times\bo$ give us the tangential projections of $-v$ and $w$, see also \cite[sec.~2]{S-Serdica}. Note that local knowledge of $\ell$ is not needed for the lens rigidity problem\footnote{If $L$ is given only, then the problem is called \textit{scattering rigidity} in some works}, and in fact, $\ell$ can be recovered locally from either $d$ or $L$, see for example the proof of  Theorem~\ref{thm_stab_global}. 

Vargo \cite{V} proved that real-analytic manifolds satisfying an additional mild condition are lens rigid.
Croke has shown that if a manifold is lens rigid, a finite quotient of it
is also lens rigid \cite{Croke04}.  He has also shown that the torus is lens
rigid \cite{Croke_scatteringrigidity}. G. Uhlmann and P.~Stefanov have shown lens rigidity locally
near a generic class of non-simple manifolds \cite{SU-lens}. In a recent work, Guillarmou [12] proved that in two dimensions, one can determine from the lens relation the conformal class of the metric if the trapped set is hyperbolic and there are no conjugate points. He also proved  deformation lens rigidity in higher dimensions under the same assumptions.
The only result we know for the lens rigidity problem with incomplete (but not local) data
is for real-analytic metric and metric close to them satisfying the micolocal condition in the next sentence \cite{SU-lens}. While in \cite{SU-lens}, the lens relation is assumed to be known on a subset only, the geodesics issued from that subset cover the whole manifold and their conormal bundle is required to cover $N^*M$. In contrast, in this paper, we have localized information. 

We state below an immediate corollary of our main result for this problem. 
For the partial data problem instead of assuming  $d=\tilde d$ locally, we can assume that $L=\tilde L$ in a neighborhood of $S_p\bo$. To reduce this problem to Theorem~\ref{thm_1} directly, we need to assume first that $c=\tilde c$ on $\bo$ near $p$ to make the definition of $\partial_\pm SM$ independent of the choice of the speed but in fact, one can redefine the lens relation in a way to remove that assumption, see \cite{SU-lens}.

\begin{theorem}\label{corollary_1}
Let $M$, $c$, $\tilde c$ be as in Theorem~\ref{thm_1} with $c=\tilde c$ on $\bo$ near $p$. Let $L=\tilde L$ near $S_p\bo$. Then $c=\tilde c$ in $M$ near $p$. 
\end{theorem}

\begin{remark}
The theorem or its corollary does not preclude the existence of an infinite set of speeds $c_j$ all  having the same boundary distance function (or lens data) in $U\times U$, where $U\subset\bo$ is some fixed small set but  not coinciding in any fixed neighborhood of $p$. In principle, this may  happen when  the maximal neighborhood of $U$, which can be covered with  strictly convex surfaces, which continuously deform $U$, shrinks when $j\to\infty$. Then the theorem does not imply existence of a fixed neighborhood of $p$, where all speeds are equal. 
If one assumes that a priori, the sound speeds have uniformly bounded derivatives of some finite order   near $p$, this situation does not arise, and this case is covered by  Theorem~\ref{thm_stab_global} below. 
\end{remark}

The linearization of the boundary rigidity and lens rigidity problem is the  tensor tomography problem, i.e., recovery of a tensor field up to  ``potential fields'' from  integrals along geodesics joining points on the boundary. It has been extensively studied in the literature for both simple and non-simple manifolds
\cite{Mu77, Dairbekov, PSU1, PSU2, PSU3, PestovSh,ShSkU, Vladimir_97, Sh1,   SU-Duke, SU-AJM, SU-caustics, Bao-Zhang}. See the book \cite{Sh-book} and  \cite{PSU3} for a recent survey. The \textit{local} tensor tomography problem has been considered
in \cite{K} for functions and real-analytic metrics and in \cite{KS} for tensors of order two and real-analytic metrics. Those results can also be thought of as support theorems of Helgason type. 
The only known results for the local problem for smooth metrics and integrals of  functions is \cite{UV:Local}.

Now we use a layer stripping type argument to obtain a global result which is different from Mukhometov's for simple manifolds.
\begin{definition}\label{def_1.1}
Let $(M,g)$ be a compact Riemannian manifold with boundary. We say that
$M$ satisfies the foliation condition by strictly convex hypersurfaces if 
 $M$ is equipped with a smooth function $\rho: {M}\to[0,\infty)$ which level sets $\Sigma_t=\rho^{-1}(t)$, $t<T$
 with some $T>0$ are strictly convex viewed from $\rho^{-1}((0,t))$ for $g$,   $d\rho$ is non-zero on these level sets,  and $\Sigma_0=\partial M$ and $M\setminus\cup_{t\in[0,T)}\Sigma_t$  has empty interior. 
\end{definition}

The statement of the global result on lens rigidity is as follows:

\begin{theorem}\label{thm_2} 
Let $n=\dim M\ge3$,  let $c>0$, $\tilde c>0$ be smooth  and equal on $\bo$,  let $\bo$  be strictly convex with respect to both $g=c^{-2}g_0$ and $\tilde g = \tilde c^{-2}g_0$. Assume that   $M$ can be foliated by strictly convex hypersurfaces for $g$. Then if  $L= \tilde L$ on $ \partial_-SM$, we have  $c=\tilde c$ in $M$.
\end{theorem}

A more general foliation condition under which the theorem would still hold is formulated in \cite{SU-Carleman}, see also Definition~\ref{def_5.1} below. In particular, $\Sigma_0$ does not need to be $\partial M$ and one can have several such foliations with the property that the  closure of their union is $M$. If we can foliate only some connected neighborhood of $\bo$, we would get $c=\tilde c$ there.  
Next,  it is enough to require that $M\setminus\cup_{t\in[0,T)}\Sigma_t$ is simple (or that it is included in a simple submanifold),   see the proof of Theorem~\ref{thm_2}   and Figure~\ref{fig:local_lens_rigidity_pic2}, to prove $c=\tilde c$ in $\cup_{t\in[0,T)}\Sigma_t$ first,  and then use Mukhometov's results to complete the proof. The class of manifolds we get in this way is larger than the simple ones, and can have conjugate points.

Speeds not necessarily radial (with $g_0$ the Euclidean metric) under the condition considered by Herglotz and Wieckert and Zoeppritz satisfy  the foliation condition of the theorem, see also Section~\ref{section_6}.   Other examples of non-simple metrics that satisfy the condition are the tubular neighborhood of a closed
geodesic in negative curvature.  These have trapped geodesics.  
It follows from the  result of \cite{RS},  that manifolds with no focal points satisfy the foliation condition. It would be interesting to know whether this is also the case for simple manifolds. As it was mentioned earlier, manifolds satisfying the foliation condition are not necessarily simple.

The linearization of the non-linear problem with partial data considered in Theorem~\ref{thm_1} was considered in \cite{UV:Local}, where uniqueness and stability were shown.  This corresponds to integrating functions along geodesics joining points in a neighborhood of $p$.
The method of proof of Theorem~\ref{thm_1} relies on using an identity proven in \cite{SU-MRL} to reduce the problem to a "pseudo-linear" one: to show uniqueness when one integrates the function $f=c^2-\tilde c^2 $ and its derivatives
on the geodesics for the metric $g$ joining points near $p$, with weight 
depending non-linearly on both $g$ and $\tilde g$. Notice that this is not a proof by linearization, and unlike the problem with full data, an attempt to do such a proof by linearization is connected with essential difficulties. 
The  proof of uniqueness for this linear transform follows the method of \cite{UV:Local} introducing an
artificial boundary and using Melrose' scattering calculus. In section~\ref{sec_2}, we do the reduction to a ``pseudo-linear problem", and in section~\ref{sec_3}, we show uniqueness for
the "pseudo-linear" problem. In section~\ref{sec_4}, we finish the proofs of the main theorems.

We also prove  H\"older conditional stability estimates related to the uniqueness theorems above.  In case of data on the whole boundary, such an estimate was proved in  \cite[section~7]{SU-JAMS} for simple manifolds and metrics not necessarily conformal to each other.   Below, the $C^k$ norm is defined in a fixed coordinate system. The next theorem is a local stability result, corresponding to the local uniqueness result in Theorem~\ref{thm_1}. 

\begin{theorem}\label{thm_stability}
There exists $k>0$ and $0<\mu<1$ with the following property. For any  $0<c_0\in C^k(M)$, $p\in\bo$, and $A>0$, there exists $\eps_0>0$ and $C>0$ with the property that for any two positive $c$, $\tilde c$ with 
\be{est}
\|c-c_0\|_{C^2} +\|\tilde c-c_0\|_{C^2} \le \eps_0, \quad \text{and} \quad \|c\|_{C^k}+ \|\tilde c\|_{C^k}\le A, 
\ee
and for any neighborhood $\Gamma$ of $p$ on $\bo$, 
we have the stability estimate
\be{stab}
\|c-\tilde c\|_{C^2(U)}\le C\|d-\tilde d\|_{C(\Gamma\times\Gamma)}^\mu
\ee
for some neighborhood $U$ of $p$ in $M$.  
\end{theorem}

In   Theorem~\ref{thm_stab_global},  we prove a H\"older conditional stability estimates  of global type as well, which can be considered as a ``stable version'' of Theorem~\ref{thm_2}.

 The plan of the paper is as follows. The reduction to a pseudo-linear problem is done in Section~\ref{sec_2}. In Section~\ref{sec_3}, we present  linear analysis using the scattering calculus. The main result there is Proposition~\ref{pr_3.2}, which is of its own interest as well. The proofs of the three uniqueness theorems are in Section~\ref{sec_4}. In Section~\ref{sec_5}, we prove the local stability result in Theorem~\ref{thm_stability} and the global Theorem~\ref{thm_stab_global}.  As an application of our results, we revisit the class of speeds studied by  Herglotz \cite{Herglotz} and Wieckert \& Zoeppritz \cite{WZ}  in Section~\ref{section_6} without assuming that they are radial and we prove that they are lens rigid. In particular, we show that their condition \r{HWZ} is equivalent to the requirement that the Euclidean spheres are strictly convex for the metric $c^{-2}\d x^2$; therefore, it is a foliation condition.  

\textbf{Acknowledgments.}  We would like to thank Christopher Croke for pointing out an error in the formulation of Theorem~\ref{thm_2} in an earlier version of the paper and for helpful comments. 

\section{Reducing the non-linear problem to a pseudo-linear one} \label{sec_2}

We recall the known fact \cite{LassasSU} that one can determine the jet of $c$ at any boundary point $p$ at which $\bo$ is convex (not necessarily strictly) from the distance function $d$ known near $(p,p)$. For a more general result not requiring convexity, see \cite{SU-lens}. Since the result in \cite{LassasSU} is formulated for general metrics, and the reconstruction of the jet is in boundary normal coordinates, we repeat the proof in this (simpler) situation of recovery a conformal factor. As in  \cite{LassasSU}, we say that $\bo$ is convex near $p\in\bo$, if for any two distinct points $p_1,p_2\in\bo$, close enough to $p$,  there  exists a geodesic $\gamma: [0,1]\to M$ joining them such that its length is $d(p_1,p_2)$ and all the interior of $\gamma$ belongs to the interior of $M$. Of course, strict convexity (positive second fundamental form at $p$) implies convexity. 

\begin{lemma}
Let $c$ and $\tilde c$ be smooth and let $\bo$ be convex at $p$ with respect to $g$ and $\tilde g$. Let $d=\tilde d$ near $(p,p)$. 
Then $\partial^\alpha c=\partial^\alpha\tilde c$ on $\bo$ near $p$ for any multiindex $\alpha$. 
\end{lemma} 

\begin{proof} 
Let $V\subset \bo$ be a neighborhood of $p$ on $\bo$ such that for any $p_1,p_2\in V$, we have the property guaranteeing convexity at $p$. Let $x^n$ be a boundary normal coordinate related to $g$; i.e., $x^n(q)=\dist(q,\bo)$, and $x^n\ge0$ in $M$. We can complete $x^n$ to a local coordinate system $(x',x^n)$, where $x'$ parameterizes  $\bo$ near $p$. 

It is enough to prove  
\be{bo}
\partial^k_{x^n} c= \partial^k_{x^n} \tilde c\quad \text{in $V$, $k=0,1,\dots$}. 
\ee
For $k=0$, this follows easily by  taking the limit in $d(p,q) = \tilde d(p,q)$, $\bo\ni q\to p$; and this can be done for any $p\in V$.  Let $l$ be the first value of $k$ for which \r{bo} fails. Without loss of generality, we may assume that it fails at $p$, and  $\partial^l_{x^n}(c-\tilde c)>0$ at $p$. Then $\partial^k_{x^n}(c-\tilde c)=0$ in $V$, $k=0,\dots,l-1$. Consider the Taylor expansion of $c-\tilde c$ w.r.t.\ $x^n$ with $x'$ close enough to $x'(p)$. We get $c-\tilde c>0$ in some neighborhood of $p$ in $M$ minus the boundary. 

Now, let $\gamma(p,q)$ be a minimizing geodesic in the metric $g$ connecting $p$ and $q$  when $q\in\bo$ as well, close enough to $p$, see also \cite{LassasSU}. Let $If(p,q)$ be the geodesic ray transform of the tensor field  $f$ defined as an  integral of $f_{ij}\dot\gamma^i\dot\gamma^j$ along $\gamma(p,q)$. All geodesics here are parameterized by a parameter in $[0,1]$ rather than being unit speed, and therefore the transform $I$ is parameterized differently than usual one. 
Then $I(g-\tilde g)>0$ by what we proved above. On the other hand, 
\[
0<I(g-\tilde g)= d^2(p,q) -I\tilde g\le d^2(p,q)  - \tilde I\tilde g = d^2(p,q)  -\tilde  d^2(p,q)=0
\]
 because $\tilde \gamma(p,q)$ minimizes  integrals of $g$ along curves connecting $p$ and $q$. This a contradiction. 
\end{proof}

The starting point is an identity in \cite{SU-MRL}. We will repeat the proof.

Let $V$, $\tilde V$ be two vector fields on a  manifold  $M$ (which will be replaced later with $S^*M$). Denote by $X(s,X\zero)$ the solution of $\dot X=V(X)$, $X(0)=X\zero$, and we use the same notation for $\tilde V$ with the  
corresponding solution are denoted by $\tilde X$. 
Then we have the following simple statement.

\begin{lemma}\label{SU-identity}
For any $t>0$ and any  initial condition $X\zero$, if $\tilde X\!\left(\cdot,X\zero\right)$ and $X\!\left(\cdot,X\zero\right) $ exist on the interval $[0,t]$, then 
\be{L1}
\tilde X\!\left(t,X\zero\right)-X\!\left(t,X\zero\right) = \int_0^t \frac{\partial \tilde X}{\partial X\zero}\!\left(t-s,X(s,X\zero)\right)\left(\tilde V-V \right)\!\left( X(s,X\zero)\right)\,ds.
\ee
\end{lemma} 
\begin{proof}
Set 
\[
F(s) = \tilde X\!\left(t-s,X(s,X\zero)\right).
\]
Then
\[
F'(s) = -\tilde V\!\left(\tilde X(t-s,X(s,X\zero))\right) + \frac{\partial \tilde X}{\partial X\zero} \!\left(t-s,X(s,X\zero)\right)V\!\left(X(s,X\zero)\right).
\]
The proof of the lemma would be complete by the fundamental theorem of calculus
\[
F(t)-F(0)=\int_0^t F'(s)\, ds
\]
 if we show the following
\begin{equation} \label{110}
 \tilde V\!\left(\tilde X(t-s,X(s,X\zero))\right)
= \frac{\partial \tilde X}{\partial X\zero}\!\left(t-s,X(s,X\zero)\right)
\tilde V\!\left(X(s,X\zero)\right).
\end{equation}
Indeed, (\ref{110}) follows from
$$
0=\left.\frac{d}{d\tau}\right|_{\tau=0} X(T-\tau,X(\tau,Z))
= -V(X(T,Z)) +\frac{\partial X}{\partial X\zero}(T,Z)
V(Z), \quad \forall T,
$$
after setting $T=t-s$, $Z= X(s,X\zero)$.
\end{proof} 

Let $c$, $\tilde c$ be two speeds. Then the corresponding metrics are $g= c^{-2} dx^2$, and  $\tilde g = \tilde c^{-2} dx^2$. The corresponding Hamiltonians and Hamiltonian vector fields are  
\[
H = \frac12 c^2g_0^{ij}\xi_i\xi_j , \qquad 
V = \left(c^2g_0^{-1} \xi, -\frac12\partial_x \left(c^2|\xi|_{g_0}^2\right)\right),
\]
and the same ones related to $\tilde c$. We used the notation  $|\xi|_{g_0}^2:= g_0^{ij}\xi_i\xi_j $.

We change the notation at this point. We denote points in the phase space $T^*M$, in a fixed coordinate system, by $z = (x,\xi)$. We denote the bicharacteristic with initial point $z$ by $Z(t,z) = (X(t,z),\Xi(t,z))$. 

Then we get the identity already used  in \cite{SU-MRL}
\be{1}
\tilde Z(t,z) - Z(t,z)= 
\int_0^t \frac{\partial\tilde Z}{\partial z} (t-s,Z 
(s,z))\big(\tilde V - V\big)(Z (s,z))\,ds.
\ee

 We can naturally think of the scattering relation $L$ and the travel time $\ell$ as functions on the cotangent bundle instead of the tangent one. Then we get the following.

\begin{proposition}\label{pr1} 
Assume  
\be{5a}
L(x_0,\xi^0) = \tilde L(x_0,\xi^0), \quad \ell(x_0,\xi^0) =\tilde \ell(x_0,\xi^0)
\ee
for some $z_0= (x_0,\xi^0)\in \partial_-S^*M$. Then 
\[ 
\int_0^{ \ell(z_0) }  \frac{\partial\tilde Z}{\partial z} ( \ell(z_0)-s,Z 
(s,z_0) )\big(V -\tilde V\big)(Z (s,z_0))\,ds =0.
\]
\end{proposition}


\subsection{Linearization near $c=1$ and $g$ Euclidean} As a simple exercise, let $c=1$, $g_{ij}=\delta_{ij}$ and linearize for $\tilde c$ near $1$ first under the assumption that $\tilde c=1$ outside an open region $\Omega\subset \R^n$. Then
\begin{equation}  \label{110'}
Z(s,z) = \mat1s01 z,\quad \frac{\partial Z(s,z)}{\partial z} = 
\mat1s01,
\end{equation}
and  we get the following formal linearization  of \r{1} 
\begin{equation} \label{112}
\int \left(f\xi-\frac12(t-s)(\partial_x f),\,
-\frac12\partial_x f \right)(x+s\xi,\xi)\,\d s=0,
\end{equation}
where 
\be{f}
f:= c^2-\tilde c^2.
\ee 
Notice that we would get the same thing if we replace $\partial \tilde Z/\partial z$ in \r{1} by $\partial  Z/\partial z$. 
We integrate over the whole line $s\in\R$ because the integrand vanishes outside the interval $[0, \ell(x,\xi)]$. 
The last $n$ components of  (\ref{112}) imply 
\begin{equation} \label{112'}
\int  \partial_x  f(x+s\xi) \,\d s=0.
\end{equation} 
Now, assume that this holds for all $(x,\xi)$. Then  $\partial_x f=0$, and since $f=0$ on $\bo$, we get $f=0$.  


\subsection{The general case} \label{sec_2.2}
We take the second $n$-dimensional component on \r{1} and use the fact that $c^2|\xi|_{g_0}^2=1$ on the bicharacteristics related to $c$. We assume that both geodesics extend to $t\in [0,\ell(x,\xi)]$.  We want to emphasize that the bicharacteristics on the energy level $H=1/2$, related to $c$, do not necessarily stay on the same energy level for the Hamiltonian $\tilde H$. We get
\be{2}
\begin{split}
&\int  \frac{\partial \tilde \Xi }{\partial x}( \ell(z)-s,Z
(s,z))(fg_0^{-1}\xi)(Z(s,z))\,\d s\\
&-\frac12  \int  \frac{\partial \tilde \Xi}{\partial \xi}( \ell(z)-s,Z 
(s,z ))(  \partial_x  (f g_0^{-1})\xi\cdot\xi)   (Z(s,z))\,\d s=0
\end{split}
\ee
for any $z\in \partial_-SM$ for which \r{5a} holds. 
As before, we integrate over  $s\in\R$ because the support of the integrand vanishes for $s\not\in [0, \ell(x,\xi)]$ (for that, we   extend the bicharacteristics formally outside so that they do not come back
). Write
\[
\partial_x  (f g_0^{-1})= g_0^{-1}\partial_x  f + (\partial_x  g_0^{-1}) f 
\]
to get
\be{2a}
\begin{split}
&\int  \frac{\partial \tilde \Xi }{\partial x}( \ell(z)-s,Z
(s,z))(fg_0^{-1}\xi)(Z(s,z))\,\d s\\
&-\frac12  \int  \frac{\partial \tilde \Xi}{\partial \xi}( \ell(z)-s,Z 
(s,z ))\left(\left(g_0^{-1}\partial_x  f+ f (\partial_x   g_0^{-1}\right)\right) \xi\cdot\xi )   (Z(s,z))\,\d s=0. 
\end{split}
\ee
One of terms on the r.h.s.\ above involves $g_0^{-1}\xi\cdot\xi$ which equals $c^{-2}$ on the bicharacteristics of $H$ on the level $1/2$.

Introduce the exit times $\tau(x,\xi)$ defined as the minimal (and the only) $t> 0$ so that $X(t,x,\xi)\in\bo$. They are well defined near $S_p\bo$, if $\bo$ is strictly convex at $p$. 
We need to write $\frac{\partial \tilde Z}{\partial z}( \ell(z)-s,Z(s,z))$ as a function of $(x,\xi)=Z(s,z)$. We have
\[
\frac{\partial \tilde Z}{\partial z}(\ell(z)-s,Z(s,z)) = \frac{\partial \tilde Z}{\partial z}(\tau(x,\xi),(x,\xi)). 
\]
Then we get, with $f$ as in \r{f},
\be{3}
J_if(\gamma):= \int \left( A_i^j(X(t),\Xi(t))(\partial_{x^j}f)(X(t))+ B_i (X(t),\Xi(t))f(X(t)) \right)\d t=0
\ee
for any bicharacteristic  $\gamma = (X(t),\Xi(t))$ (related to the speed $c$) in our set, where
\be{4}
\begin{split}
A_i^j\left(x,\xi\right) = &-\frac12 \frac{\partial\tilde  \Xi_i}{\partial \xi_j}(\tau(x,\xi),(x,\xi)) c^{-2}(x),
\\
B_i\left(x,\xi\right) = &\frac{\partial \tilde\Xi_i}{\partial x^j}(\tau(x,\xi),(x,\xi))g_0^{ik}(x) \xi_k  -\frac12 \frac{\partial\tilde  \Xi_i}{\partial \xi_j}(\tau(x,\xi),(x,\xi))  (\partial_{x^j}   g_0^{-1}(x))\xi\cdot\xi.
\end{split} 
\ee

A major inconvenience with this representation is that the exit time function $\tau(x,\xi)$ (recall that we assume  strong convexity) becomes singular at $(x,\xi)\in T^*\bo$. More precisely, the normal derivative w.r.t.\ $x$ when $\xi$ is tangent to $\bo$ has a square root type of singularity. On the other hand, we have some freedom to extend the flow outside $M$ since we know that the jets of $c$ and $\tilde c$ at $\bo$ are the same near $p$: therefore, any smooth local extension of $c$ is also a smooth extension of $\tilde c$. Then for any $(x,\xi)\in \partial_-S^*M$ close enough to $S^*_{x_0}M$, the bicharacteristics originating from it will be identical once they exit $T^*M$ but are still close enough to it. Similarly, instead of starting from $T^*\bo$, we could start at points and codirections close to it, but away from $\bar M$. 

\begin{figure}[h!] 
  \centering
  \includegraphics[scale=0.85]{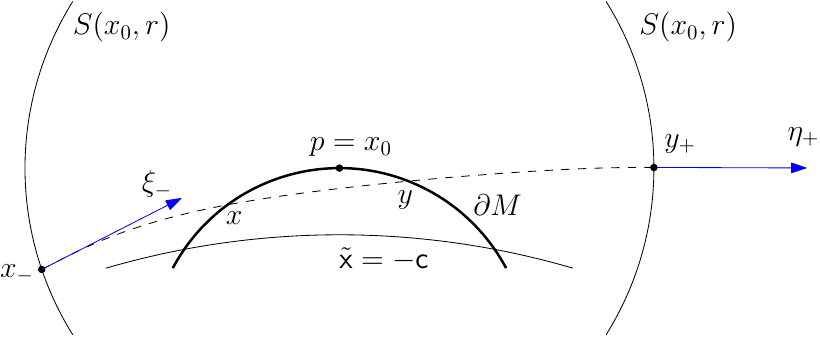}
\caption{The redefined scattering relation $(x_-,\xi_-)\mapsto (y_+,\eta_+)$}
  \label{fig:local_lens_rigidity_pic1}
\end{figure}

With this in mind, we  push the boundary  away a bit. Let $x_0$ represent the point $p$ near which we work, in a fixed coordinate system. Extend $g_0$ smoothly near $x_0$. 
Let $S(x_0,r)$ be the sphere in the metric $c^{-2} \d x^2$ centered at $x_0$ with radius $0<r\ll1$. For $(x,\xi)$ with $x$ in the geodesic ball $B(x_0,r)$, redefine $\tau(x,\xi)$ to be the travel time from $(x,\xi)$ to $S(x_0,r)$. Let $U_-\subset\partial_-SB(x_0,r)$ be the set of all points on  $S(x_0,r)$ and incoming unit directions so that the corresponding geodesic in the metric $g$ is close enough to one tangent to $\bo$ at $x_0$. Similarly, let $U_+$ be the set of such pairs with outgoing  directions. 
Redefine the scattering relation $L$ locally to act from $U_-$ to $U_+$, and redefine $\ell$ similarly, see Figure~\ref{fig:local_lens_rigidity_pic1}.  
Then under the assumptions of Theorem~\ref{corollary_1}, 
$L=\tilde L$ and $\ell=\tilde \ell$ on $U_-$. We can apply the construction above by replacing $\partial_\pm SM$ locally by $U_\pm$. The advantage we have now is that on $U_-$, the travel time $\tau$ is non-singular. Equalities \r{3}, \r{4} are preserved then. 

We now have
\be{5}
\begin{split}
A_i^j(x_0,\xi) &= -\frac12 \frac{\partial\tilde  \Xi_i}{\partial \xi_j}(r,(x_0,\xi)) c^{-2}(x_0),\\
 B_i(x_0,\xi) &= \frac{\partial \tilde\Xi_i}{\partial x^j}(r,(x_0,\xi)) g_0^{ik}(x_0) \xi_k   -\frac12 \frac{\partial\tilde  \Xi_i}{\partial \xi_j}(r,(x,\xi))  (\partial_{x^j}   g_0^{-1}(x_0))\xi\cdot\xi.
\end{split} 
\ee
Then by the strict convexity, 
\be{7}
A_i^j(x,\xi) =  -\frac12 c^{-2}(x)\delta_i^j +O(r), \quad \text{for $(x,\xi)\in S^*M $ near $S^*_p\bo$}.
\ee

\subsection{A new linear problem} 
The arguments above lead to the following linear problem: 

\textbf{Problem.} 
Assume \r{3} holds with some $f$ supported in $M$, for all geodesics  close to the ones originating from $S^*_{x_0}\partial M$ (i.e. initial point $x_0$ and all unit initial co-directions tangent to $\bo$). Assume that $\bo$ is strictly convex at $x_0$ w.r.t.\ the speed $c$. Assume \r{7}.  
 Is it true that $f=0$ near $x_0$? 

We show below in Proposition~\ref{pr_3.2}, that the answer is affirmative. Note that this reduces the original non-linear problem to a linear one but this is not a linearization. Then  Theorem~\ref{corollary_1} follows from it. On the other hand, Theorem~\ref{corollary_1} is not equivalent to that problem because the weight there has a specific structure, thus Proposition~\ref{pr_3.2} is a more general statement.

\section{Linear analysis}\label{sec_3}
We first recall the setting introduced in \cite{UV:Local} in our
current notation. There the scalar X-ray transform along geodesics was
considered, namely for $\beta\in SM$,
$$
(If)(\beta)=\int_{\RR}f(\gamma_\beta(t))\,\d t,
$$
where $\gamma_\beta$ is the geodesic with lift to $SM$ having
starting point $\beta\in SM$. Here $M$ is assumed to have a strictly
convex boundary, which can be phrased as the statement that if $\rho$
is a defining function for $\pa M$, then
$-\frac{\d^2}{\d t^2}(\rho\circ\gamma_\beta)|_{t=0}\geq C>0$ whenever
$\beta\in S\pa M$. One then considers a point $p\in \pa M$, and
another function $\tilde\foliation$, denoted in \cite{UV:Local} by $\tilde
x$, such that $\tilde\foliation(p)=0$, $\d\tilde\foliation(p)=-\d\rho(p)$, and the
level sets of $\tilde\foliation$ near the $0$ are strictly concave when
viewed from the superlevel sets (which are on the side of $p$ when
talking about the $\level$-level set with $\level<0$), i.e.\
$\frac{\d^2}{\d t^2}(\tilde\foliation\circ\gamma_\beta)|_{t=0}\geq C>0$ if $\beta\in
S\{\tilde\foliation=\level\}$, namely if
$\frac{\d}{\d t}(\tilde\foliation\circ\gamma_\beta)|_{t=0}=0$.
For $\level>0$, we denote $M\cap\{\tilde\foliation>-\level\}$ by $\Omega_\level$; we assume
that $\level_0>0$ is such that
$\overline{\Omega_{\level_0}}$ is compact on $M$, and the concavity assumption
holds on $\Omega_{\level_0}$. Then it was shown that the X-ray
transform $I$ restricted to $\beta\in S\overline{\Omega_{\level}}$
such that $\gamma_\beta$ leaves $\Omega_{\level}$ with both endpoints
on $\pa M$ (i.e.\ at $\rho=0$) is injective if $\level>0$ is sufficiently small, and indeed one has a stability
estimate for $f$ in terms of $If$ on exponentially weighted
spaces.

To explain this in detail, let
$\foliation=\foliation_\level=\tilde\foliation+\level$ be the boundary
defining function of the artificial boundary,
$\tilde\foliation=-\level$, that we introduced; indeed it is
convenient to work in $\tilde M$, a $\CI$ manifold extending $M$
across the boundary, extending $\tilde\foliation$ to $\tilde M$ smoothly, and defining
$\Omegaext=\Omegaext_\level=\{\foliation_\level>0\}$
as the extension of $\Omega$, so $\Omegaextb$ is a smooth manifold
with boundary, with only the artificial boundary being a boundary
face. Then one writes
$\beta=(\lambda,\omega)=\lambda\pa_{\foliation}+\omega\pa_{\loccoord}\in
S_{\foliation,\loccoord}\tilde M$ relative to a product
decomposition $(-\level_0,\level_0)_\foliation\times U_\loccoord$ of
$\tilde M$ near $p$.
The concavity condition becomes that for $\beta$ whose
$\lambda$-component vanishes,
$$
2\alpha(\foliation,\loccoord,0,\omega)=\frac{\d^2}{\d t^2}(\foliation\circ\gamma_\beta)|_{t=0}\geq 2C>0,
$$
with a new $C>0$,
see the discussion preceding Equation~(3.1) in \cite{UV:Local}.
For $\chi\in\CI_0(\RR)$, $\chi\geq 0$,
$\chi(0)>0$, one considers the map
$$
L_0v(\foliation,\loccoord)=\int_\RR\int_{\sphere^{n-2}} \foliation^{-2}\chi(\lambda/\foliation)v(\foliation,\loccoord,\lambda,\omega)\,\d\lambda\,\d\omega
$$
defined for $v$ a function on $S_{\Omegaextb}\tilde M$. This differs
from \cite{UV:Local}  in that the weight $\foliation^{-2}$ differs by $1$
from the weight $\foliation^{-1}$ used in \cite{UV:Local}; this simply has the effect of removing an
$\foliation^{-1}$ in \cite[Proposition~3.3]{UV:Local}, as compared to
the proposition stated below. If $\level$ is sufficiently
small, or instead $\chi$ has sufficiently small support, for
$(\foliation,\loccoord)\in\overline{\Omega}$,
$I$ only integrates over
points in $\beta\in S\Omegaextb$
such that $\gamma_\beta$ leaves $\Omega_{\level}$ with both endpoints
on $\pa M$, i.e.\ over $\beta$ corresponding to {\em
  $\Omega_\level$-local geodesics} -- the set of the latter is denoted
by $\cM_\level$. We refer to the discussion around
\cite[Equation~(3.1)]{UV:Local} for more detail, but roughly speaking
the concavity of the level sets of $\foliation$ means that the
geodesics that are close to being tangent to the foliation, with
`close' measured by the distance from the artificial boundary,
$\foliation=0$, then they cannot reach $\foliation=0$ (or reach again, in case
they start there) without reaching $\foliation=\level'$ for some fixed
$\level'>0$; notice that the geodesics involved in the integration
through a point on the level set $\foliation=\tilde\level$
make
an angle $\lesssim\tilde\level$ with the tangent space of the level
set due to the compact support of $\chi$.
Then we consider the map
$P=L_0\circ I$. The main technical result of \cite{UV:Local},
whose notation involving the so-called scattering Sobolev spaces
$\Hsc^{s,r}(\Omegaextb)$ and scattering pseudodifferential operators
$\Psisc^{s,r}(\Omegaextb)$ is explained below, was:

\begin{proposition}(See \cite[Proposition~3.3 and Lemma~3.6]{UV:Local}.)\label{prop:scalar-psdo}
For $\digamma>0$ let
$$
P_\digamma=e^{-\digamma/x} P e^{\digamma/x}:\CI_0(\Omegaextb)\to\CI(\Omegaextb).
$$
Then $P_\digamma\in\Psisc^{-1,0}(\Omegaextb)$.

Further, if $\level>0$  is sufficiently small, then for
suitable choice of $\chi\in\CI_0(\RR)$ with $\chi(0)=1$, $\chi\geq 0$,
$P$ is elliptic
in $\Psisc^{-1,0}(\Omegaextb)$ on a neighborhood of
$\overline{\Omega}$.

Shrinking $\level$ further if needed, $P_\digamma$ satisfies the
estimate
\begin{equation}\label{eq:scalar-lin-est}
\|v\|_{\Hsc^{s,r}(\Omegaextb)}\leq C\|P_\digamma v\|_{\Hsc^{s+1,r}(\Omegaextb)}
\end{equation}
for $v$ supported in $M\cap\Omegaextb$.
\end{proposition}

We now briefly explain the role of the so-called scattering
pseudodifferential operators and the corresponding Sobolev spaces
(which are typically used to study phenomena `at infinity') in
our problem (where there is no obvious `infinity'); we refer to
\cite[Section~2]{UV:Local} for a more thorough exposition. These concepts were introduced by Melrose, \cite{RBMSpec},
in a general geometric setting, but on $\RR^n$ these operators actually
correspond to a special case of H\"ormander's Weyl calculus \cite{Hor}, also
studied earlier by Shubin \cite{Sh} and Parenti \cite{Pa}. So consider
the reciprocal spherical coordinate map, $(0,\ep)_\foliation\times
\sphere^{n-1}_\theta\to\RR^n$, with
$\sphere^{n-1}\subset\RR^n$ (unit sphere), and map
$(\foliation,\theta)\mapsto \foliation^{-1}\theta\in\RR^n$. This
map is a diffeomorphism onto its range, and it provides a
compactification of $\RR^n$ (the so called radial or geodesic compactification) by adding
$\{0\}_\foliation\times \sphere^{n-1}_\theta$ as infinity to $\RR^n$,
to obtain $\overline{\RR^n}$, which is now diffeomorphic to a
ball. Now for general $U_\loccoord$ above, we may regard, at least
locally\footnote{That is, possibly at the cost of
shrinking it; in fact all concepts below are essentially local within
$\sphere^{n-1}$, thus even in full generality one can reduce
scattering objects to (conic regions near infinity in) $\RR^n$ this way, much as standard Sobolev
spaces and pseudodifferential operators are so reduceable to subsets
of $\RR^n$ with compact closure} $U_\loccoord$ also as a coordinate chart in
$\sphere^{n-1}$, and thus obtain an identification of
$\overline{\Omega}$ with a region intersection $\overline{\RR^n}$,
thus our artificial boundary $\foliation=0$ corresponds to infinity at $\RR^n$.
In particular, notions from $\RR^n$ can now be transferred to a
neighborhood of our artificial boundary. Since the relevant vector
fields on $\RR^n$ are generated by translation invariant vector
fields, which are complete under the exponential map, the transferred
analysis replaces the incomplete geometry of standard vector fields on
$\Omegaextb$ by a complete one. Concretely, these vector
fields, when transferred, become linear combinations of
$\foliation^2\pa_\foliation$ and $\foliation\pa_{\loccoord_j}$, with
smooth coefficients. In particular, these are the vector fields with
respect to which Sobolev regularity is measured. Thus,
$\Hsc^{s,r}(\Omegaextb)$ is the so-called {\em scattering Sobolev
space}, which is locally, under the above identification, just the
standard weighted Sobolev space $H^{s,r}(\RR^n)$, see
\cite[Section~2]{UV:Local}, while $\Psisc^{s,r}(\Omegaextb)$ is
Melrose's {\em scattering pseudodifferential algebra}, which locally, again
under this identification, simply
corresponds to quantizing symbols $a$ with $|D_z^\alpha D_\zeta^\beta
a|\leq C\langle z\rangle^{r-|\alpha|}\langle\zeta\rangle^{s-|\beta|}$
on $T^*\RR^n=\RR^n_z\times\RR^n_\zeta$,
see again \cite[Section~2]{UV:Local} for more detail. Note that
ellipticity in this algebra, called full ellipticity, is {\em both}
in the sense as $|z|\to\infty$ and $|\zeta|\to\infty$, i.e.\ modulo
symbols one order lower {\em and} with an extra order of decay as $|z|\to\infty$.

Notice that \eqref{eq:scalar-lin-est} implies the estimate
$$
\|f\|_{e^{\digamma/\foliation} \Hsc^{s,r}(\Omegaextb)}\leq C\|P f\|_{e^{\digamma/\foliation}\Hsc^{s+1,r}(\Omegaextb)}
$$
for the unconjugated operator, valid when $f$ is supported in $M\cap\Omegaextb$. Rewriting $P$ as $L_0\circ I$, this
gives that for $\delta>0$, $s\geq -1$,
$$
\|f\|_{e^{(\digamma+\delta)/\foliation} H^{s}(\Omegaextb)}\leq C\|If\|_{e^{\digamma/\foliation} H^{s+1}(\cM_\level)};
$$
see the discussion in \cite{UV:Local} after Lemma~3.6.

After this recollection,
we continue by generalizing \eqref{3} to regard the functions
$\pa_{x_j}f$ and $f$ entering into it as independent unknowns, while
restricting the transform to the region of interest $\Omega=\Omega_\level$. So let $\tilde J_i$ be defined by
$$
\tilde J_i (u_0,u_1,\ldots,u_n)(\beta):= \int_{\gamma_\beta} \left( A_i^j(X(t),\Xi(t))u_j(X(t))+ B_i (X(t),\Xi(t))u_0(X(t)) \right)\d t,
$$
where $\gamma_\beta$ is the geodesic with lift to $S \Omega$ having
starting point $\beta\in S \Omega$. Let $\tilde J=(\tilde J_1,\ldots,\tilde J_n)$.
This is a vector valued version of the geodesic X-ray transform
considered in \cite{UV:Local}, and described above, sending functions on $\Omega$ with values in
$\Cx^{n+1}$ to functions with values in $\Cx^n$.
We then define $L$ as a map from
$\Cx^n$-valued functions on $S\Omega$ to $\Cx^n$ valued functions on $\Omega$ by
$$
Lv(\foliation,\loccoord)=\int_\RR\int_{\sphere^{n-2}} \foliation^{-2}\chi(\lambda/\foliation)v(\foliation,\loccoord,\lambda,\omega)\,\d\lambda\,\d\omega
$$
as in \cite{UV:Local}; this is a diagonal operator: $L=L_0\otimes\Id$. Then we consider the map
$P=L\circ\tilde J$, and in addition to the properties mentioned above
in the scalar setting, we are also interested in continuity properties
in terms of the background data, such as in the background metric as
well as the function $\foliation$.
Recall that the
map
\begin{equation}\label{eq:Gamma+-def}
\Gamma_+:S\tilde M\times[0,\infty)\to[\tilde M\times\tilde M;\diag],\
\Gamma_+(\foliation,\loccoord,\lambda,\omega,t)=((\foliation,\loccoord),\gamma_{\foliation,\loccoord,\lambda,\omega}(t))
\end{equation}
is a local diffeomorphism, and similarly for $\Gamma_-$ in which
$(-\infty,0]$ takes the place of $[0,\infty)$; see the discussion
around \cite[Equation~(3.2)-(3.3)]{UV:Local}; indeed this is true for
more general curve families.  Here $[\tilde
M\times\tilde M;\diag]$ is the {\em blow-up} of $\tilde M$ at the
diagonal $z=z'$, which essentially means the introduction of spherical/polar
coordinates, or often more conveniently projective coordinates, about
it. Concretely, writing the (local) coordinates from the two factors of
$\tilde M$ as $(z,z')$,
\begin{equation}\label{eq:blowup-coords}
z,|z-z'|,\frac{z-z'}{|z-z'|}
\end{equation}
give (local) coordinates on this space. Note that when $\Gamma_\pm$ are given by
geodesics of a metric $g$ just $C^1$-near a fixed background metric $g_0$, as $\CI$ maps,
$\Gamma_\pm$ depend continuously on $g$ in the $\CI$ topology.

In order to consider continuity properties in $\tilde\foliation$ in a
$C^1$-neighborhood of a fixed function $\tilde\foliation_0$, it is
convenient to use the map $(\tilde\foliation,\loccoord)$ to identify a
neighborhood of $p$ with a neighborhood $O$ of the origin in $\RR\times
U$. Thus, for $\level$ fixed, but $\tilde\foliation$ being $C^1$-close to
$\tilde\foliation_0$, on this fixed background $O\subset \RR\times U$, the pulled pack metrics $\underline{g}$
depend continuously on $(\tilde\foliation,g)$ as maps
$\CI(\tilde M)\times\CI(\tilde M;S^2\tilde M)\to\CI(O)$; this normalizes $\tilde\foliation$
to be simply the first coordinate function $\tufoliation$ on $O$. Correspondingly,
below, the continuous dependence, of all objects discussed, on
$\tilde\foliation$ (in the $\CI$ topology) is automatic: what is meant
always is that by pull-back to $O$ the resulting objects, living on
fixed domains such as $\{\tufoliation+\level>0\}$,
depend continuously on $g$ and $\tilde\foliation$, which follows from
the continuous dependence of these objects on $\underline{g}$. Since we do not want to
overburden the notation, we do {\em not} write this pull-back explicitly.

The main technical result here is:

\begin{proposition}\label{prop:elliptic-op}
For $\digamma>0$,  let
$$
P_\digamma=e^{-\digamma/\foliation} P e^{\digamma/\foliation}:\CI_0(\Omegaextb;\Cx^{n+1})\to\CI(\Omegaextb;\Cx^{n}).
$$
Then $P_\digamma\in\Psisc^{-1,0}(\Omegaextb;\Cx^{n+1},\Cx^n)$, and the
map $\Gamma_\pm\mapsto P_\digamma$ is continuous from the $\CI$
topology to the Fr\'echet topology of $\Psisc^{-1,0}(\Omegaextb;\Cx^{n+1},\Cx^n)$.

Further, if $\level$ is sufficiently small and \r{7} holds,  then for
suitable choice of $\chi\in\CI_0(\RR)$ with $\chi(0)=1$, $\chi\geq 0$,
if we write
$P_\digamma=(P_0,\tilde P)$, with $P_0$ corresponding to the first component,
$\tilde P$
the last $n$ components, in the domain space, then $\tilde P$ is elliptic
in $\Psisc^{-1,0}(\Omegaextb;\Cx^n,\Cx^n)$ in a neighborhood of $\overline{\Omega}$.
\end{proposition}

\begin{remark}
Notice that ellipticity being an open condition, this means that there
exists $\level_0>0$ such that if $\level<\level_0$, then the same $\chi$
works for all $g$ $\CI$-close to $g_0$.

Further, in view of the paragraph preceding the proposition, the map
$(\tilde\foliation,\Gamma_\pm)\mapsto P_\digamma$ is continuous from the $\CI$
topology to the Fr\'echet topology of
$\Psisc^{-1,0}(\Omegaextb;\Cx^{n+1},\Cx^n)$, where the latter is
understood to actually stand for
$\Psisc^{-1,0}(\{\tufoliation+\level\geq 0\};\Cx^{n+1},\Cx^n)$ via the
identifications discussed above.
\end{remark}

\begin{proof}
This is simply a vector valued version of
\cite[Proposition~3.3]{UV:Local} and \cite[Lemma~3.6]{UV:Local},
recalled above in Proposition~\ref{prop:scalar-psdo}. In
particular, to show $P\in\Psisc^{-1,0}(\Omegaextb;\Cx^{n+1},\Cx^n)$, it
suffices to show that $P$ is a matrix of pseudodifferential operators
$P_{ij}\in\Psisc^{-1,0}(\Omegaextb)$, $i=1,2,\ldots,n$,
$j=0,1,2,\ldots,n$, depending continously on $\Gamma_\pm$. But for $j>0$, with $j=0$ being completely
analogous, $P_{ij}w$ has the form
$$
\int_\RR\int_{\sphere^{n-2}} \foliation^{-2}\chi(\lambda/\foliation )\int A_i^j(X_{\foliation,\loccoord,\lambda,\omega}(t),\Xi_{\foliation,\loccoord,\lambda,\omega}(t))w(X_{\foliation,\loccoord,\lambda,\omega}(t))\,\d t\,\,\d\lambda\,\d\omega
$$
The only difference from \cite[Proposition~3.3]{UV:Local} then is the
presence of the weight factor
$$A_i^j(X_{\foliation,\loccoord,\lambda,\omega}(t),\Xi_{\foliation,\loccoord,\lambda,\omega}(t)).$$ 
It is convenient to rewrite this via the metric identification, say by
$g_0$, in terms of tangent vectors. Changing the notation for the
flow, in our coordinates $(\foliation,\loccoord,\lambda,\omega)$,
writing now
$$
(\gamma_{\foliation,\loccoord,\lambda,\omega}(t),\gamma'_{\foliation,\loccoord,\lambda,\omega}(t))
=(\Foliation_{\foliation,\loccoord,\lambda,\omega}(t),\Loccoord_{\foliation,\loccoord,\lambda,\omega}(t),\Lambda_{\foliation,\loccoord,\lambda,\omega}(t),\Omega_{\foliation,\loccoord,\lambda,\omega}(t))
$$
for
the lifted geodesic $\gamma_{\foliation,\loccoord,\lambda,\omega}(t)$,
$$
\tilde
A_i^j(\Foliation_{\foliation,\loccoord,\lambda,\omega}(t),\Loccoord_{\foliation,\loccoord,\lambda,\omega}(t),\Lambda_{\foliation,\loccoord,\lambda,\omega}(t),\Omega_{\foliation,\loccoord,\lambda,\omega}(t))
$$
replaces
$A_i^j(X_{\foliation,\loccoord,\lambda,\omega}(t),\Xi_{\foliation,\loccoord,\lambda,\omega}(t)).$
As in
\cite[Proposition~3.3]{UV:Local} one rewrites the integral in terms of
coordinates $(\foliation,\loccoord,\foliation',\loccoord')$ on the left and right factors of $\Omegaextb$
(i.e.\ one explicitly expresses the Schwartz kernel), using that the
map $\Gamma_+$ of \eqref{eq:Gamma+-def}
is a local diffeomorphism, and similarly for $\Gamma_-$; we again
refer to the discussion
around \cite[Equation~(3.2)-(3.3)]{UV:Local}.
Further,
$$
(\foliation,\loccoord,\lambda,\omega,t)\mapsto\gamma'_{\foliation,\loccoord,\lambda,\omega}(t)
$$
is a smooth map
$S M\times\RR\to\RR^n$ (depending continuously on $\Gamma_\pm$ in the
respective $\CI$ topologies) so composing it with $\Gamma_\pm^{-1}$
from the right, one can re-express the integral giving
$P_{ij}w$ away from
the boundary as
$$
\int w(z')|z-z'|^{-n+1}b\Big(z,\frac{z-z'}{|z-z'|},|z-z'|\Big)\,\d z'
$$
as in \cite[Equation~(3.7)]{UV:Local}, with $b$ a smooth function of
the indicated variables (thus smooth on $[\tilde
M\times\tilde M;\diag]$), depending continuously on $\Gamma_\pm$ in
the respective $\CI$ topologies, and with
$$
b\Big(z,\frac{z-z'}{|z-z'|},0\Big)=\tilde\chi\Big(z,\frac{z'-z}{|z'-z|}\Big)\tilde
A_i^j\Big(z,\frac{z'-z}{|z'-z|}\Big)\sigma\Big(z,\frac{z'-z}{|z'-z|}\Big)
$$
with $\sigma>0$, bounded below by a positive constant, a weight factor,
and where $\chi(\lambda/\foliation)$ is written as
$\tilde\chi(\foliation,\loccoord,\lambda,\omega)$.
Recall from \cite[Section~2]{UV:Local} that coordinates on Melrose's
scattering double space, on which the Schwartz kernels of elements of
$\Psisc^{s,r}(\Omegaextb)$ are conormal to the diagonal, near the
lifted scattering diagonal, are (with $\foliation\geq 0$)
$$
\foliation,\ \loccoord,\
X=\frac{\foliation-\foliation'}{\foliation^2},\ Y=\frac{\loccoord-\loccoord'}{\foliation}.
$$
Further, it is convenient to write coordinates on $[\tilde M\times
\tilde M;\diag]$ in the region of interest (see the beginning of the paragraph of
Equation~(3.10) in \cite{UV:Local}), namely (the lift of)
$|\foliation-\foliation'|<C|\loccoord-\loccoord'|$, as
$$
\foliation,\loccoord,|\loccoord-\loccoord'|,\frac{\foliation-\foliation'}{|\loccoord-\loccoord'|},\frac{\loccoord-\loccoord'}{|\loccoord-\loccoord'|},
$$
with the norms being Euclidean norms,\footnote{This
  is an example of partial projective coordinates for a blow-up.}
instead of \eqref{eq:blowup-coords}; we write $\Gamma_\pm$ in terms of
these. Note that these are $\foliation,\loccoord,\foliation|Y|,\frac{\foliation|X|}{|Y|},\hat
Y$.
Then, similarly, near the boundary as
in \cite[Equation~(3.13)]{UV:Local}, one obtains the Schwartz kernel
\begin{equation}\begin{aligned}\label{eq:SK-near-ff}
K^\flat(\foliation,\loccoord,X,Y)=\sum_\pm e^{-\digamma
  X/(1+\foliation X)}&\chi\Big(\frac{X}{|Y|}+|Y|\tilde\Lambda_\pm\Big(\foliation,\loccoord,\foliation|Y|,\frac{\foliation|X|}{|Y|},\hat
Y\Big)\Big)\\&
\tilde A_i^j\Big(\Gamma_\pm^{-1}\Big(\foliation,\loccoord,\foliation|Y|,\frac{\foliation X}{|Y|},\hat Y\Big)\Big)|Y|^{-n+1}J_\pm\Big(\foliation,\loccoord,\frac{X}{|Y|},|Y|,\hat Y\Big),
\end{aligned}\end{equation}
with the density factor $J$ smooth, positive, depending continuously
on $\Gamma_\pm$ in the respective $\CI$ topologies, $=1$ at $\foliation=0$.
Here
$$
\foliation,\loccoord,|Y|,\frac{X}{|Y|},\hat Y
$$
are valid
coordinates on the blow-up of the scattering diagonal in\footnote{This
  is another example of partial projective coordinates for a blow-up.} $|Y|>\ep|X|$,
$\ep>0$, which is the case automatically on the support of the kernel
due to the argument of $\chi$, cf.\ the discussion after \cite[Equation(3.12)]{UV:Local}, so the argument of
$\tilde A_i^j\circ\Gamma_\pm^{-1}$ is smooth {\em on this blown up
  space} still depending continuously on $\Gamma_\pm$ in the
respective $\CI$ topologies. We can evaluate this
argument: for instance, by
\cite[Equation(3.10)]{UV:Local},
$$
(\Lambda\circ\Gamma_\pm^{-1})\Big(\foliation,\loccoord,\foliation|Y|,\frac{\foliation
  X}{|Y|},\hat Y\Big)=
\foliation\frac{X}{|Y|}+\foliation|Y|\tilde\Lambda_\pm\Big(\foliation,\loccoord,\foliation|Y|,\frac{\foliation
  X}{|Y|},\hat Y\Big)
$$
with $\tilde\Lambda_\pm$ smooth, while the subsequent equation in the same
location gives
$$
(\Omega\circ\Gamma_\pm^{-1}) \Big(\foliation,\loccoord,\foliation|Y|,\frac{\foliation
  X}{|Y|},\hat Y\Big)=\hat Y+\foliation|Y|\tilde\Omega_\pm\Big(\foliation,\loccoord,\foliation|Y|,\frac{\foliation
  X}{|Y|},\hat Y\Big)
$$
with $\tilde\Omega_\pm$ smooth; here both $\tilde\Lambda_\pm$ and
$\tilde\Omega_\pm$ depend continuously on $\Gamma_\pm$ in the respective $\CI$ topologies.
This proves the first part of the
proposition as in \cite[Proposition~3.3]{UV:Local}.

To prove the second part, note that in view of \eqref{eq:SK-near-ff}
(which just needs to be evaluated at $\foliation=0$), \cite[Lemma~3.5]{UV:Local} is
  replaced by the statement that the boundary principal symbol of
  $P_{ij}$ in
  $\Psisc^{-1,0}(\Omegaextb)$ is twice the $(X,Y)$-Fourier transform
of
\begin{equation}\label{eq:Pij-symbol}
e^{-\digamma
  X}\chi\Big(\frac{X-\alpha(0,\loccoord,0,\hat Y)|Y|^2}{|Y|}\Big)
\tilde A_i^j\Big(0,\loccoord,0,0,\hat Y\Big)|Y|^{-n+1},
\end{equation}
while for $P_{i0}$ it is twice the $(X,Y)$-Fourier transform of
$$
e^{-\digamma
  X}\chi\Big(\frac{X-\alpha(0,\loccoord,0,\hat Y)|Y|^2}{|Y|}\Big)
\tilde B_i\Big(0,\loccoord,0,0,\hat Y\Big)|Y|^{-n+1},
$$
with $\tilde B_i$ defined analogously to $\tilde A_i^j$.
(Recall that $2\alpha(\foliation,\loccoord,\lambda,\omega)$ is the $\foliation$ component of
$\gamma''_{\foliation,\loccoord,\lambda,\omega}(0)$, and the convexity assumption on $\foliation$
is that $\alpha$ is positive; see \cite{UV:Local} above Equation~(3.1).)
For $A_i^j=-\frac{1}{2}c^{-2}(x_0)\delta_i^j$, see \r{7}, the invertibility of
the principal symbol, with values in $n\times n$ matrices, of the
principal symbol of $\tilde P$ follows when $\chi$ is chosen as in \cite[Lemma~3.6]{UV:Local}, for it is
$-\frac{1}{2}c^{-2}(x_0)$ times the boundary symbol in
\cite[Lemma~3.6]{UV:Local} times the $n\times n$ identity matrix. In
general, due to the perturbation stability of the property of
invertibility, the same follows for $\level$ sufficiently small.  
\end{proof}

\begin{corollary}\label{cor:elliptic-est}
With the notation of Proposition~\ref{prop:elliptic-op}, there is $\tilde\level>0$
such that if $0<\level<\tilde\level$, then $P_\digamma$ satisfies the
estimate
\begin{equation}\label{eq:vector-lin-est}
\|u\|_{\Hsc^{s,r}(\Omegaextb)}\leq C\|P_\digamma u\|_{\Hsc^{s+1,r}(\Omegaextb)}+C\|u_0\|_{\Hsc^{s,r}(\Omegaextb)}
\end{equation}
for $u$ supported in $M\cap\Omegaextb$, with the constant $C$ uniform
in $\level$. Further, fixing $s,r$, there exist $k$ and $\ep>0$ and $\rho_0<0$ such that
$C$ can be taken
uniform for $\Gamma_\pm$ $\ep$-close to a reference $\Gamma_\pm^0$ in
$C^k$, and the estimate holds even for $u$ supported in $\rho\geq \rho_0$
(in place of $\rho\geq 0$).
\end{corollary}

\begin{remark}
As in the case of the preceding proposition, the dependence on
$\tilde\foliation$ is also continuous, i.e.\ by possibly increasing $k$, $C$
can be taken uniform for $\Gamma_\pm$ $C^k$-close to a reference
$\Gamma_\pm^0$ and $\tilde\foliation$ $C^k$-close to a reference $\tilde\foliation_0$.
\end{remark}

\begin{proof}
By the density of elements of $\dCI(M\cap\Omegaextb)$ in
$\Hsc^{s,r}(\Omegaextb)$ supported in $M\cap\Omegaextb$, it suffices
to consider $u\in \dCI(M\cap\Omegaextb)$ to prove \eqref{eq:vector-lin-est}.

Consider $s=0$, $r=0$.
Let $\Lambda\in\Psisc^{1,0}(\Omegaextb)$ be elliptic and invertible
(one can e.g.\ locally identify $\Omega$ with $\RR^n$; then on the Fourier transform side
multiplication by $\langle\xi\rangle$ works). Thus, $\Lambda
P_\digamma$, with $\Lambda$ acting diagonally on $\dCI(\Omegaextb;\Cx^n)$, is in
$\Psisc^{0,0}(\Omegaextb;\Cx^{n+1},\Cx^n)$, depending continuously, in the Fr\'echet topology of pseudodifferential
operators, on
$\Gamma_\pm$ in $\CI$, and $\Lambda\tilde
P\in \Psisc^{0,0}(\Omegaextb;\Cx^{n},\Cx^n)$ is elliptic in
$\Psisc^{0,0}(\Omegaextb;\Cx^{n},\Cx^n)$ {\em locally in a
  neighborhood of $\Omega$}. This implies, as presented in
\cite{UV:Local} after Lemma~3.6 (without the uniform discussion on $\Gamma_\pm$), relying on the arguments at the end of Section~2
there, that there exist $k$ and $\ep>0$ such that if $\Gamma_\pm$ is
$\ep$-close to $\Gamma_\pm^0$ in $C^k$, then if $\level>0$ is sufficiently
small then $\Lambda\tilde
P$ satisfies
\begin{equation}\label{eq:local-elliptic-bound}
\|w\|_{L^2_\scl(\Omegaextb)}\leq C_0\|\Lambda\tilde P w\|_{L^2_\scl(\Omegaextb)}
\end{equation}
for $w$ supported in $\overline{\Omega}$. Here
$L^2_\scl(\Omegaextb)=\Hsc^{0,0}(\Omegaextb)$ is the $L^2$ space
relative to a non-degenerate {\em scattering
  density} ---  the latter are equivalent to the lifted Lebesgue measure from
$\RR^n$, thus are bounded
multiples of $\frac{\d\foliation\,\d\loccoord}{\foliation^{n+1}}$.

We recall the essential part of this argument briefly. One considers
the whole family of domains $\Omegaextb_{\level}$, which can be
identified with each other locally in the region of interest by the
maps
$\Phi_{\level}(\tilde\foliation,\loccoord)=(\tilde\foliation+\level,
\loccoord)$, i.e.\ simply translation in the
$\tilde\foliation$-coordinate, so instead of considering a family of
spaces with an operator on each of them, one can consider a fixed
space, denote this by $\Omegaextb_0$, with a continuous (in $\level$) family of operators,
$T_\level$, namely $T_\level=(\Phi_\level^{-1})^*\Lambda\tilde P
\Phi_\level^*$; these depend continuously on $\Gamma_\pm$ in the
Fr\'echet sense discussed above. Notice that we are interested in the
region\footnote{If we are interested in the region $\rho\geq\rho_0$,
  $\rho_0<0$, within $\Omegaextb$, with $\rho_0$ to be specified, we
  take e.g.\ $\rho_0=-\level$ in the discussion below, and the desired continuous
  $f$ with $f(0)=0$, and in this case with $\foliation_\level\leq
  f(\level)$ on $\Omegaext_\level\cap\{\rho\geq\rho_0\}$, still exists.}
$\Omega_\level$, and that there is a continuous function $f$ on $\RR$
with $f(0)=0$ and $\foliation_\level\leq f(\level)$ on
$\Omega_\level$. Correspondingly, in the translated space,
$\tilde\foliation\leq f(\level)$ on the image of $\Omega_\level$;
notice that this region shrinks as $\level>0$ goes to $0$.
On the other hand, there is a fixed open set $O\subset\Omegaextb$, a
neighborhood of $x_0$, on
which the operators $T_\level$ are elliptic in
$\Psisc^{0,0}(\Omegaextb)$ for $0\leq |\level|<\level_0$. Let $K_0$ be
a compact subset of $O$, still including a neighborhood of $x_0$,
$\phi\in\CI_0(O)$ be identically $1$ on a neighborhood of $K$.
Then the elliptic parametrix construction (which is local, and uniform
in $\level$ by the continuity) produces a parametrix family
$G_\level\in\Psisc^{0,0}(\Omegaextb;\Cx^{n},\Cx^n)$, $G_\level
T_\level=\Id+E_\level$, where $E_\level\in\Psisc^{0,0}(\Omegaextb;\Cx^n,\Cx^n)$
only, but $\phi E_\level\phi\Psisc^{-\infty,-\infty}(\Omegaextb;\Cx^n,\Cx^n)$, uniformly in
$\level$, depending continuously on $\Gamma_\pm$ in the Fr\'echet topologies. Then (multiplying the parametrix identity by $\phi$ from
left and right and applying to $v$) for $v$ supported in $K$, $(\Id+\phi
E_\level\phi)v=\phi G_\level T_\level v$. Now, the Schwartz kernel of
$\phi E_\level \phi$ is Schwartz, i.e.\ is bounded by
$C_N(\foliation\foliation')^N$ for any $N$, uniformly in $\level$. (Here
we write, say, the Schwartz kernel relative to scattering
  densities, but
as $N$ is arbitrary, this makes little difference.) For $\level>0$, let
$\phi_\level$ be supported, say, in $\tilde\foliation\leq 2f(\level)$,
identically $1$ near the region $\tilde\foliation\leq f(\level)$; one
may assume that $\phi\equiv 1$ on $\supp\phi_\level$ by making
$\level>0$ small.
Then, by
Schur's lemma, $\phi_\level\phi E_\level \phi\phi_\level$, acting say on
$L^2_\scl(\Omegaextb)$ (i.e.\ the $L^2$-space relative to scattering
densities) is bounded by $C'_N f(\level)^N$ for any $N$. Thus, there is
$\level_1>0$ such that
$\Id+\phi_\level\phi E_\level \phi\phi_\level=\Id+\phi_\level E_\level\phi_\level$ is invertible for
$0<\level<\level_1$ on $L^2_\scl(\Omegaextb)$. (Notice here that this
requires the smallness of a seminorm of $\phi_\level\phi E_\level
\phi\phi_\level$ in the Fr\'echet topology of pseudodifferential
operators; the continuity discussed above means that that this
requires the $C^k$-closeness of $\Gamma_\pm$ to $\Gamma_\pm^0$ for
some $k$.) In particular,
for $v$ supported in $\tilde\foliation\leq f(\level)$, so $\phi_\level v=v$,
$(\Id+\phi_\level
E_\level\phi_\level)v=\phi_\level G_\level T_\level v$, so inverting the
factor on the left and then undoing the
transformation $\Phi_\level$ gives the desired conclusion \eqref{eq:local-elliptic-bound}.

Thus, with $u=(u_0,\tilde u)$, we have, with all norms being
$L^2_\scl(\Omegaextb)$-norms,
$$
\|\Lambda P_\digamma u\|^2=\|\Lambda P_0 u_0+\Lambda \tilde P\tilde
u\|^2
=\|\Lambda\tilde P\tilde u\|^2+\langle\Lambda\tilde
P \tilde u,\Lambda P_0 u_0\rangle+\langle \Lambda P_0 u_0,\Lambda\tilde
P \tilde u\rangle+\langle\Lambda P_0 u_0,\Lambda P_0 u_0\rangle.
$$
By \eqref{eq:local-elliptic-bound},
$\|\Lambda\tilde
P \tilde u\|\geq C_0\|\tilde u\|$, $C_0>0$. On the other hand,
$\|\Lambda P_0 u_0\|\leq C_1\|u_0\|$, $\|\Lambda \tilde P \tilde
u\|\leq C_1\|\tilde u\|$ as elements of
$\Psisc^{0,0}(\Omega)$ are $L^2$-bounded. (By the continuity of the
map from $\Gamma_\pm$ to $\Lambda P_0$ in the appropriate Fr\'echet
topologies, again there is $k$ such that $C_1$ is uniform when
$\Gamma_\pm$ is $C^k$-close to $\Gamma_\pm^0$.) Using the Cauchy-Schwartz
inequality, for $\delta>0$,
$$
|\langle \Lambda P_0u_0,\Lambda\tilde
P \tilde u\rangle|\leq \frac{\delta}{2} \|\tilde
u\|^2+\frac{C_1^2}{2\delta}\|u_0\|^2.
$$
Thus,
the last three terms are bounded by $\delta\|\tilde u\|^2+(1+\delta^{-1})C_1^2\|u_0\|^2$ in
absolute value, so we conclude that, with $\delta=C_0^2/2$,
\begin{equation}\label{eq:tilde-u-est}
\frac{C_0^2}{2}\|\tilde u\|^2\leq \|\Lambda P_\digamma
u\|^2+C_1^2(1+2C_0^{-2})\|u_0\|^2,
\end{equation}
completing the proof of the corollary if $s=r=0$.

The general case follows via conjugating $P_\digamma$ by an elliptic,
invertible, element of $\Psisc^{s,r}(\Omegaextb)$, which is thus an
isomorphism from $\Hsc^{s,r}(\Omegaextb)$ to
$L^2=\Hsc^{0,0}(\Omegaextb)$. Note that such a conjugation does not
change the principal symbol, thus the ellipticity.
\end{proof}

We now remark that the even simpler setting of the {\em scalar transform
with a positive weight $A$ on $S^*M$}, $I_A$, which was not considered in
\cite{UV:Local}, and which can
be considered a special case of $\tilde J$ with $(u_0,u_1,\ldots,u_n)$
simply replaced by $u$. Thus, for consistency with the above notation,
let $\tilde A(\foliation,\loccoord,\lambda,\omega)$ be the weight on
$SM$ induced by the metric identification, and
let
$$
(I_Af)(\beta)=\int_{\RR}\tilde A(\gamma_\beta(t),\gamma'_\beta(t))f(\gamma_\beta(t))\,\d t,\
\beta\in SM.
$$
Then the above argument gives that
$P_\digamma\in\Psisc^{-1,0}(\Omegaextb)$. If at $S_{\pa M}M$,  the
weight $A$ is independent of the momentum variable $\xi$, it further
gives that $P_\digamma$ is elliptic in a neighborhood of $\overline{\Omega}$.
More generally, a modification of the argument of \cite{UV:Local} due to
H.~Zhou \cite{HZhou:Local} allows one to show that the principal
symbol is fully elliptic in $\overline{\Omega}$ in the scattering sense merely assuming that $A$ is
positive (but {\em not} the independence condition just mentioned). To see this,
one has to Fourier transform \eqref{eq:Pij-symbol} in $(X,Y)$ with
$\tilde A_i^j$ replaced by $\tilde A$. The $X$-Fourier transform is
unaffected by the presence of $\tilde A$, and gives, as in \cite[Equation~(3.16)]{UV:Local},
$$
|Y|^{2-n}e^{-\alpha(\digamma+i\xi)|Y|^2}\hat\chi(\xi-i\digamma)\tilde
A(0,\loccoord,0,\hat Y).
$$
Replacing $\chi$ with a Gaussian, $\chi(s)=e^{-s^2/(2\nu)}$,
$\nu=\digamma^{-1}\alpha$, which does not have compact support, but an
approximation argument (in symbols of order $-1$) will give this desired property, one can compute that this is, up to a
constant factor,
$$
\digamma^{-1/2}\alpha^{1/2}|Y|^{2-n}e^{-\digamma^{-1}(\xi^2+\digamma^2)\alpha|Y|^2/2}\tilde
A(0,\loccoord,0,\hat Y).
$$
We need to compute the Fourier transform in $Y$. Following
\cite{HZhou:Local}, one expresses this in polar coordinates in $Y$:
\begin{equation*}\begin{aligned}
&\digamma^{-1/2}\int_{\sphere^{n-2}}\int_0^\infty e^{-i|Y|\hat Y\cdot\eta}\alpha^{1/2}e^{-\digamma^{-1}(\xi^2+\digamma^2)\alpha|Y|^2/2}\tilde
A(0,\loccoord,0,\hat Y)\,\d|Y|\,\d\hat
Y\\
&=\frac{1}{2}\digamma^{-1/2}\int_{\sphere^{n-2}}\int_{\RR}e^{-it\hat
  Y\cdot\eta}\alpha^{1/2}e^{-\digamma^{-1}(\xi^2+\digamma^2)\alpha t^2/2}\tilde
A(0,\loccoord,0,\hat Y)\,\d t\,\d\hat Y,
\end{aligned}\end{equation*}
which in turn becomes, up to a constant factor,
$$
(\xi^2+\digamma^2)^{-1/2}\int_{\sphere^{n-2}}e^{-\digamma|\hat Y\cdot\eta|^2/(2(\xi^2+\digamma^2)\alpha)}\tilde
A(0,\loccoord,0,\hat Y)\,\d\hat Y.
$$
The integrand is now positive, which gives the desired ellipticity at $\foliation=0$. (One
also needs to check the ellipticity as $(\xi,\eta)\to\infty$; this is
standard, see \cite{UV:Local,HZhou:Local}.) One
proceeds with an approximation argument as in
\cite{UV:Local,HZhou:Local} to complete the proof of the ellipticity.
Thus, the above argument gives the estimate
$$
\|f\|_{\Hsc^{s,r}(\Omegaextb)}\leq C\|P_\digamma f\|_{\Hsc^{s+1,r}(\Omegaextb)}.
$$
As discussed in \cite{UV:Local} after Lemma~3.6, this yields the
following corollary:

\begin{corollary}
The weighted scalar transform with a positive weight $A$ on $S^*M$,
with $\tilde A$ the associated weight on $SM$,
$$
(I_Af)(\beta)=\int_{\RR}\tilde A(\gamma_\beta(t),\gamma'_\beta(t))f(\gamma_\beta(t))\,\d t,
$$
satisfies that
for $\digamma>0$ there is $\level_0>0$ such that for $\delta>0$,
$0<\level<\level_0$, $s\geq -1$, we have
$$
\|f\|_{e^{(\digamma+\delta)/\foliation} H^{s}(\Omega_\level)}\leq C\|I_Af\|_{e^{\digamma/\foliation} H^{s+1}(\cM_\level)}.
$$
\end{corollary}

We now return to the actual case of interest and apply Corollary~\ref{cor:elliptic-est}
with $u_0=e^{-\digamma/\foliation} f$,
$u_j=e^{-\digamma/\foliation}\pa_j f$. If we show that (uniformly in
$\Gamma_\pm$, which are indeed irrelevant for this argument) given $\tilde\delta>0$ there is $\level_0>0$ such that $\|u_0\|_{L^2}\leq
\tilde\delta\|\tilde u\|_{L^2}$, i.e.\ $e^{-\digamma/\foliation} f$ is bounded
by a small multiple of a derivative of $f$ times
$e^{-\digamma/\foliation}$ in $L^2$, when $f$ is supported in
$\overline{\Omega_\level}$, $0<\level<\level_0$,
then for $\tilde\delta>0$ sufficiently small \eqref{eq:vector-lin-est}
proves that if
$P_\digamma u$ vanishes, then so does $u$, i.e.\ in this case so does $f$, for
the $u_0$ term can then be absorbed into the left hand side of
\eqref{eq:vector-lin-est}:
\begin{equation}\label{eq:vector-lin-est-mod}
(1-\tilde\delta C)\|u\|_{\Hsc^{0,0}(\Omegaextb)}\leq C\|P_\digamma u\|_{\Hsc^{1,0}(\Omegaextb)}.
\end{equation}
Further, rewriting this by removing the weights
$e^{-\digamma/\foliation}$, and estimating the norms in terms if the
standard $L^2$-based space, cf.\ the discussion after
\cite[Lemma~3.6]{UV:Local} already referenced above,\footnote{The
  $\delta$ loss is actually just the loss of a power of $\foliation$,
  due to change of the measure.} gives, for $\delta>0$,
\begin{equation}\label{eq:f-transform-est-0}
\|f\|_{e^{(\digamma+\delta)/\foliation} H^1(\Omega_\level)}\leq
C\|\tilde J(f,\pa_1 f,\ldots,\pa_n f)\|_{e^{\digamma/\foliation} H^{1}(\cM_\level)},
\end{equation}
with $C$ uniform for $\Gamma_\pm$ $C^k$-close to $\Gamma_\pm^0$.

But
this can now be easily done: let $V$ be a smooth vector field with
$V\foliation=0$, so $V$ is tangent
to the boundary of $\Omega_\level$ for every $\level$, and make the
non-degeneracy assumption that, for some $\level_1>0$,
there is a continuous $T:[0,\level_1)\to\RR$ such that $T(0)=0$ and
the $V$-flow takes every point in $\overline{\Omega_{\level}}$ to $\pa M$
in time $\leq T(\level)$ (i.e. outside
the original manifold). Then the Poincar\'e inequality for $V$
gives\footnote{Notice that our treatment of the
  $\tilde\foliation$-dependence of the problem relies on reducing to a
  model, where $\tilde\foliation$ is replaced by a fixed function $\tufoliation$, so
  the following argument is in fact uniform in $\tilde\foliation$.}
\begin{equation}\label{eq:Poincare}
\|v\|_{L^2_\scl}\leq C_2 T(\level)\|Vv\|_{L^2_\scl},
\end{equation}
for $v$ vanishing outside $M$, hence the constant is
small if $T(\level)>0$ is small. (Here the $L^2$ space we need is the
scattering $L^2$-space, $L^2_\scl$, which
is $\foliation^{(n+1)/2}$ times the standard $L^2$-space, but the extra weight
does not affect the argument, since $V$ commutes with multiplication
by powers of $\foliation$.)

To see \eqref{eq:Poincare}, we recall a standard proof of
the local Poincar\'e inequality: in order to reduce confusion with the notation, let
$(z_1,...,z_n)=(z',z_n)$ be the coordinates, $z_1=0$ being the boundary of
$\Omega$ (so $\foliation$ would be $z_1$), and assume that the flow of
$\pa_{z_n}$ flows from every point in $\Omega$ to outside
the region in `time' $\leq \delta$. To normalize
the argument, assume that $z_n\geq 0$ in $\Omega$, and we want to
estimate $v$ in
$z_n\leq \delta$. We assume that the $L^2$ space is given by a density
$F(z')\,|dz|$. Then, for $v\in\CI(\RR^n_+)$, $\RR^n_+=[0,\infty)_{z_1}\times\RR^{n-1}$, with support in $z_n\geq
0$, by the fundamental theorem of calculus and the Cauchy-Schwartz inequality,
\begin{equation*}\begin{aligned}
|v(z',z_n)|=\Big|\int_0^{z_n} \pa_n v(z',t)\,\d t\Big|
&\leq\Big(\int_0^{z_n} 1\,\d t\Big)^{1/2}\Big(\int_0^{z_n} |\pa_n v(z',t)|^2\,\d t\Big)^{1/2}\\
&\leq \delta^{1/2}\Big(\int_0^\delta |\pa_n v(z',t)|^2\,\d t\Big)^{1/2}.
\end{aligned}\end{equation*}
Squaring both sides, multiplying by $F(z')$, and integrating in $z',z_n$ (to $z_n=\delta$) gives
$$
\int_{z_n\leq \delta} |v(z',z_n)|^2\,F(z')\,\d z
\leq \delta^2\int_{t\leq \delta} |\pa_n v(z',t)|^2\,\d t\,F(z')\,\d z'.
$$
This says that actually
$$
\|v\|_{L^2(\RR^n_+;F(z')\,|\d z|)}\leq \delta\|\pa_n v\|_{L^2(\RR^n_+;F(z')\,|\d z|)},
$$
proving the claim (using $F(z')=z_1^{-n-1}$) in view of the quasi-isometry invariance (which
gives a constant factor) of the
bound \eqref{eq:Poincare}. Even if there is more complicated topology,
so there are no global coordinates and vector fields as stated,
dividing up the problem into local pieces and adding them together
gives the desired result: taking steps of size $\delta$, one needs
$T/\delta$ steps to cover the set, using cutoff functions to localize
is easily seen to give a bound proportional to $T$.

On the other hand, in view of the
strict convexity of the boundary, one can construct such a $V$ and $T$.
With
$v=e^{-\digamma/\foliation}f$, this is exactly the desired conclusion
since $V(e^{-\digamma/\foliation}f)=e^{-\digamma/\foliation}Vf$,
namely \eqref{eq:vector-lin-est-mod}, and thus
\eqref{eq:f-transform-est-0}, hold.

In fact, we can prove the analogue of
\eqref{eq:vector-lin-est-mod} on stronger spaces:
\begin{equation}\label{eq:vector-lin-est-mod-gen}
\|u\|_{\Hsc^{s,0}(\Omegaextb)}\leq C\|P_\digamma
u\|_{\Hsc^{s+1,0}(\Omegaextb)},\qquad s\geq 0,
\end{equation}
which in turn gives, for $s\geq 0$,
\begin{equation}\label{eq:f-transform-est}
\|f\|_{e^{(\digamma+\delta)/\foliation} H^{s+1}(\Omega_\level)}\leq
C\|\tilde J(f,\pa_1 f,\ldots,\pa_n f)\|_{e^{\digamma/\foliation} H^{s+1}(\cM_\level)},
\end{equation}
with $C$ uniform for $\Gamma_\pm$ $C^k$-close to $\Gamma_\pm^0$. To
see this, notice that \eqref{eq:vector-lin-est} is an elliptic
estimate when we have $u_0=e^{-\digamma/\foliation} f$,
$u_j=e^{-\digamma/\foliation}\pa_j f$, for 
$$
\|u\|_{\Hsc^{s,r}(\Omegaextb)}\leq C\|P_\digamma u\|_{\Hsc^{s+1,r}(\Omegaextb)}+C\|u_0\|_{\Hsc^{s,r}(\Omegaextb)}
$$
implies that\footnote{Notice that the scattering derivatives are actually
weaker than the derivatives $\pa_j$ entering $u$, and one can absorb
the term given by commuting a scattering derivative through
$e^{-\digamma/\foliation}$ into the last term on the right hand side.}
$$
\|e^{-\digamma/\foliation} f\|_{\Hsc^{s+1,r}(\Omegaextb)}\leq
C\|P_\digamma
e^{-\digamma/\foliation}(f,df)\|_{\Hsc^{s+1,r}(\Omegaextb)}+C\|e^{-\digamma/\foliation
 } f\|_{\Hsc^{s,r}(\Omegaextb)}.
$$
For $s=1$, $r=0$, we have the second term on the right hand side
controlled by $P_\digamma
e^{-\digamma/\foliation}(f,df)$ in $\Hsc^{1,0}$ in view of
\eqref{eq:vector-lin-est-mod}, so the $\Hsc^{2,0}$-norm of
$e^{-\digamma/\foliation} f$ is in fact controlled by $\|P_\digamma
e^{-\digamma/\foliation}(f,df)\|_{\Hsc^{2,0}(\Omegaextb)}$, which we
can iterate further\footnote{We can easily allow $s\geq 0$
  non-integer by slightly modifying the argument here.} to prove \eqref{eq:vector-lin-est-mod-gen}.

In summary we have proved that with $L^2(\Omega)=L^2(\Omega_\level)$ the standard $L^2$-space
now (as the exponential weight $e^{-\digamma/\foliation}$ maps such $f$ to
$L^2_\scl(\Omegaextb)$, see also the discussion after Lemma~3.6 in \cite{UV:Local}):

\begin{proposition}\label{pr_3.2}
There is $\level_0>0$ such that for $0<\level<\level_0$, if
$f\in H^1(\Omega_\level)$ and $\tilde J (f,\pa_1 f,\ldots,\pa_n
f)=0$, then $f=0$.

In fact, for $\digamma>0$, $s\geq 1$, there exist $\level_0>0$,
$k$ and $\ep>0$ such that the following holds. For $\delta>0$
there is $C>0$ such that if
$0<\level<\level_0$, $\Gamma_\pm$ is $\ep$-close to
$\Gamma_\pm^0$ in $C^k$, $\tilde\foliation$ is $\ep$-close to
$\tilde\foliation_0$ in $C^k$, then
$$
\|f\|_{e^{(\digamma+\delta)/\foliation} H^s(\Omega_\level)}\leq
C\|\tilde J(f,\pa_1 f,\ldots,\pa_n f)\|_{e^{\digamma/\foliation} H^{s}(\cM_\level)}.
$$

Moreover, with $\Omega_{\level,\rho_0}=\Omegaext_\level\cap\{\rho\geq
\rho_0\}$, and $\cM_{\level,\rho_0}$ being defined analgously to
$\cM_\level$ with $\pa M=\{\rho=0\}$ being replaced by
$\{\rho=\rho_0\}$, we have: for $\digamma>0$ and $s\geq 1$ there exist $\level_0>0$, $\rho_0<0$,
$k$ and $\ep>0$ such that the following holds. For $\delta>0$
there is $C>0$ such that if
$0<\level<\level_0$, $\Gamma_\pm$ is $\ep$-close to
$\Gamma_\pm^0$ in $C^k$, $\tilde\foliation$ is $\ep$-close to
$\tilde\foliation_0$ in $C^k$, then $f\in
H^{s+1}(\Omega_{\level,\rho_0})$
implies that
$$
\|f\|_{e^{(\digamma+\delta)/\foliation} H^s(\Omega_{\level,\rho_0})}\leq
C\|\tilde J(f,\pa_1 f,\ldots,\pa_n f)\|_{e^{\digamma/\foliation} H^{s}(\cM_{\level,\rho_0})}.
$$
\end{proposition}

\section{Proof of Theorems~\ref{thm_1}-\ref{thm_2}} \label{sec_4}

\begin{proof}[Proof of Theorem~\ref{corollary_1}]
Let $c$ and $\tilde c$ be as in the theorem. Redefine the scattering relation $\mathcal{L}$ as in Figure~\ref{fig:local_lens_rigidity_pic1}. By Proposition~\ref{pr1}, we get $J_jf(\gamma)=0$, see \r{3} for all geodesics close enough to the ones tangent to $\bo$ at $p$. The weights are given by \r{5}, in the new parameterization, with the ellipticity condition  satisfied by \r{7}. Then Proposition~\ref{pr_3.2} implies $f=0$ in a neighborhood of $p$, where $f=c^2-\tilde c^2$ as in \r{f}. 
\end{proof}

\begin{proof}[Proof of Theorem \ref{thm_1}]
Note first that we can complete $(M,c^{-2}_0)$, and similarly $(M,\tilde c^{-2}g_0)$ to compact Riemannian manifolds without boundary. Then we can choose a neighborhood of $p$ small enough so that the exponential map based on any point of that neighborhood is a local diffeomorphism for short enough vectors, both for $c$ and for $\tilde c$. This implies that there is a neighborhood $U$ of $p$ so that the distance $d(p_1,p_2)$ between any two points $p_1$, $p_2$ in $U$ is realized by $|\exp_{p_1}^{-1}(p_2)|$, related to the first and the second metric, respectively, where $\exp$ is the localized exponential map as above. We can easily recover $c$ and $\tilde c$ on $\bo\cap U$ by taking the limit $p_1\to p_2$. 
As Michel proved \cite{Michel}, for simple manifolds, the scattering elation can be recovered by differentiating the distance function, see also \cite[sec.~2]{S-Serdica}. This applies to our case as well because if $U$ is small enough. Then the proof follows from Theorem~\ref{corollary_1}. 
\end{proof}

\begin{proof}[Proof of Theorem \ref{thm_2}]
Theorem \ref{thm_2} is now an easy consequence of Theorem \ref{thm_1} using a layer stripping argument. Let $f=c^2-\tilde c^2$.
Assume $f\neq 0$, then $\supp f$ has non-empty
interior. On the other hand, let $\tau=\inf_{\supp f}\rho$; if $\tau= T$ we
are done, for then $\supp f\subset M\setminus\cup_{t\in[0,T)}\Sigma_t$. Thus, suppose $\tau<T$, so
$f\equiv 0$ on $\Sigma_t$ for $t<\tau$, but there exists $x\in\Sigma_\tau\cap\supp f$ (since $\supp f$ is closed).
We will show below how to use Theorem~\ref{corollary_1} on $M_\tau:= \rho^{-1}(\tau,\infty)$  to conclude that a neighborhood of $x$ is disjoint from $\supp f$ to obtain a contradiction. 

All we need to show is that the scattering relations $L_\tau$ and $\tilde L_\tau$ on $\Sigma_\tau$ coincide. Note that $\Sigma_\tau = \bo_\tau$ is strictly convex for $\tilde g$ as well because the second fundamental form for $\tilde g$ can be computed by taking derivatives from the exterior $\rho<\tau$, where $g=\tilde g$. Fix $(x_\tau,v_\tau)\in \partial_-SM_\tau$, see Figure~\ref{fig:local_lens_rigidity_pic2}.  The geodesic $\gamma_{x_\tau,v_\tau}(s)$  cannot hit $\Sigma_\tau$ again for negative ``times'' $s$ because otherwise, we would get a contradiction with the strict convexity at $\Sigma_t$, where $t$ corresponds to the smallest value of  $\rho$ on that geodesic between two contacts with $\Sigma_\tau$. 
Since $c=\tilde c$ outside $M_\tau$,  $\gamma_{x_\tau,v_\tau}(s)$  and $\tilde \gamma_{x_\tau,v_\tau}(s)$   coincide outside $M_\tau$ for  $s<0$.   Proposition~\ref{pr_nont} below shows that this negative  geodesic ray must be non-trapping, i.e.,  $\gamma_{x_\tau,v_\tau}$ would hit $\bo$ for a finite negative time $s$  at some point and direction $(x,v)\in \partial_-SM$. In the same way, we show that the same holds for the positive part, $s>0$, of a geodesic issued from $L_\tau(x_\tau,y_\tau) =:(y_\tau,w_\tau)\in \partial_+SM_\tau$; and the corresponding point on $\partial_+SM$ will be denoted by $(y,w)$. Then, since $L(x,v)= (y,w)$,  we would also get $L_\tau(x_\tau,v_\tau) =  (y_\tau,w_\tau)= \tilde L_\tau(x_\tau,v_\tau)$. 
\begin{figure}[h!] 
  \centering
  \includegraphics[scale=0.7]{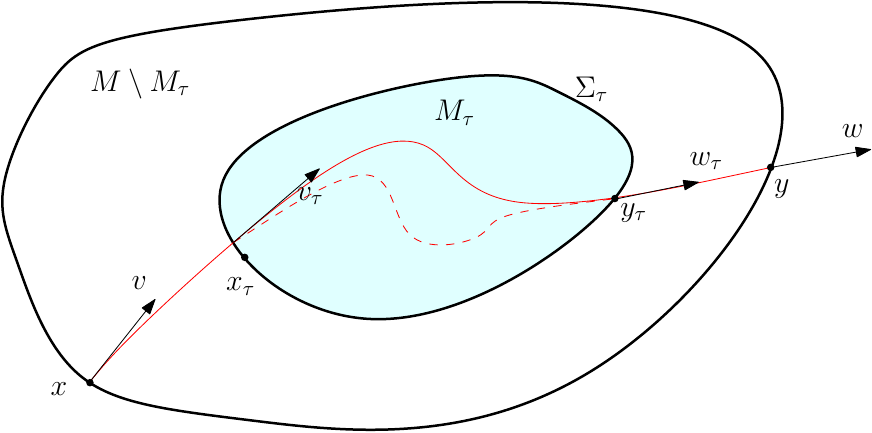}
\caption{One can recover the scattering relation on $\Sigma_\tau$  knowing that on $\bo$.}
  \label{fig:local_lens_rigidity_pic2}
\end{figure}
\end{proof}

\section{Stability Analysis}\label{sec_5}
In this section, we prove the  stability estimate in Theorem~\ref{thm_stability} and   a global stability estimate, see Theorem~\ref{thm_stab_global} below. We follow some of the ideas in  \cite[section~7]{SU-JAMS}.

\subsection{Boundary stability}

We start with stability at the boundary. 

\begin{theorem}  \label{thm_st_bd}
Let $c$ and $\tilde c$ be such that $\bo $ is strictly convex both w.r.t. \ $c^{-2}g_0$ and $\tilde c^{-2}g_0$, and $\Gamma\subset \subset \Gamma'\subset \bo$ be two sufficiently small open subsets of the boundary. Then 
\be{bdry}
\big\| \partial_{x^n}^k  
(c-  \tilde c)\big\|_{C^m(\bar\Gamma)} \le C_{k,m} \big\| d ^2 -\tilde d^2 \big\|_{C^{m+2k+2} \big(\overline{\Gamma'\times \Gamma'}\big)},
\ee
where $C_{k,m}$ depends only on $\Omega$ and on  a upper bound of $c$, $c^{-1}$, $\tilde c$, $\tilde c^{-1}$   in $C^{m+2k+5}(M)$.
\end{theorem}

\begin{proof}
We know from \cite[section~7]{SU-JAMS}  that it is true for the metrics in boundary normal coordinates. More precisely, let  $(x',x^n)$ be boundary local coordinates for $g$, i.e., locally,
\[
c^{-2}g_0 =h_{\alpha\beta}\d x^\alpha \d x^\beta+(\d x^n)^2 =: h_{ij}\d x^i\d x^j.
\]
Note that the notational convention is different than the one in section~\ref{sec_3}; $(\mathsf{x},\mathsf{y})$ is now replaced by $(x',x^n)$, and $x$ and $y$ are generic points in (the local chart) in  $\R^n$. 
 We use the convention that Greek indices run from $1$ to $n-1$, while Latin ones run from $1$ to $n$. 
Let $\psi$ be the diffeomorpishm fixing the boundary pointwise near $p$, i.e., $\psi(x',0)= (x',0)$, 
 so that
\be{a1}
\psi^*(\tilde c^{-2}g_0) =\tilde h_{\alpha\beta}\d x^\alpha \d x^\beta+(\d x^n)^2
\ee
near $p$. 
Then \r{bdry} holds for $h-\tilde h$. 
 
We have, with $y=\psi(x)$, 
\be{a2}
 \psi^*(\tilde c^{-2}g_0)(x) =(c^2 \tilde c^{-2})(y(x)) h_{i'j'}(y(x)) \frac{\partial y^{i'}}{\partial x^{i}} \frac{\partial y^{j'}}{\partial x^{j}}\d x^{i}\d x^{j}.
\ee
In particular, by \r{a1} and \r{a2},
\be{a3}
(c^2\tilde c^{-2})(y(x))h_{i'j'}(y(x)) \frac{\partial y^{i'}}{\partial x^{\alpha}} \frac{\partial y^{j'}}{\partial x^{\beta}} = \tilde h_{\alpha\beta}(x),
\ee
which can be written as
\be{a3a}
\left((c^2\tilde c^{-2})(y(x))-1\right) h_{i'j'}(y(x)) \frac{\partial y^{i'}}{\partial x^{\alpha}} \frac{\partial y^{j'}}{\partial x^{\beta}} = \tilde h_{\alpha\beta}(x) -  h_{i'j'}(y(x)) \frac{\partial y^{i'}}{\partial x^{\alpha}} \frac{\partial y^{j'}}{\partial x^{\beta}} .
\ee
We have
\[
\partial y^i/\partial x^\beta=\delta^i_\beta \quad \text{for $x^n=0$},
\]
which, of course, can be further differentiated w.r.t.\ tangential variables $x^\alpha$. 
 Therefore, 
\be{a3'}
 \left(\tilde c^{-2}c^2 -1\right)h_{\alpha\beta}|_{x^n=0} = \big(\tilde h_{\alpha\beta}- h_{\alpha\beta}\big)|_{x^n=0}.
\ee
Since $\{h_{\alpha\beta}\}$ is invertible (at lest for one multi-index is enough), we get \r{bdry} for $\tilde c^{-2}c^2-1 $, and therefore for $\tilde c-c$ for $k=0$. 

In order to do the same for $k=1$, we need to estimate $\partial (y-x)/\partial x^n$ at $x^n=0$ first. The diffeomorphism $\psi$ identifies boundary normal coordinates for $g$ and $\tilde g$. The normal $\nu=(0,1)$ (here, $0$ is $n-1$ dimensional), is unit for $g=c^{-2}g_0$ but it has length $\tilde c^{-1} c$ in the metric $\tilde g$. The inner  unit vector in the metric $\tilde g$ is therefore $\tilde \nu = \tilde c c^{-1}(0,1)$, hence
\be{a3''}
y=\tilde \gamma_{(x',0),(0,\tilde c c^{-1})}(x^n).
\ee
Differentiate this at $x^n=0$ to get that it depends only on $\tilde c c^{-1}$ on $x^n=0$ but not on its derivatives; in fact,
\be{a4}
\frac{\partial y}{\partial x^n}\bigg|_{x^n=0}= (0,\tilde c c^{-1}),
\ee
and this can be differentiated w.r.t.\ $x'$. 
Note that we can get the same result by comparing the $\alpha n$ metric elements  in \r{a1} and \r{a2}. Then $\partial (y-x)/\partial x^n|_{x^n=0}=  (0,\tilde c c^{-1}-1)$ 
By what we proved, the latter satisfies \r{bdry}. 

Differentiate \r{a3a} with respect to $x^n$ at $x^n=0$. Using \r{a4}, we get
\be{a5}
\begin{split}
\tilde c c^{-1} 
h_{\alpha\beta}  \frac{\partial (c^2\tilde c^{-2}-1)}{\partial x^n}\Big|_{x^n=0} 
=&
-\left( c^2\tilde c^{-2}( x',0) -1\right)  \frac{\partial }{\partial x^n}\Big|_{x^n=0} \left(  h_{i'j'}(y(x)) \frac{\partial y^{i'}}{\partial x^{\alpha}} \frac{\partial y^{j'}}{\partial x^{\beta}} 
 \right)\\
&\quad+\frac{\partial }{\partial x^n}\Big|_{x^n=0} \left(\tilde h_{\alpha\beta}(x) -  h_{i'j'}(y(x)) \frac{\partial y^{i'}}{\partial x^{\alpha}} \frac{\partial y^{j'}}{\partial x^{\beta}}\right).
 \end{split}
\ee
To estimate the last derivative, notice that by \r{a4},
\[
\frac{\partial^2 y}{\partial x^n \partial x^\alpha}\bigg|_{x^n=0}= \frac{\partial }{\partial x^\alpha}\bigg|_{x^n=0}(0,\tilde c c^{-1}),
\]
and we proved the estimate for the tangential derivatives of $\tilde c c^{-1}$ (same as the tangential derivatives of $\tilde c c^{-1}-1$) already. Therefore, the second summand on the r.h.s.\ of \r{a5} can be written as a sum of a term involving the normal derivative of $h-\tilde h$, which estimate is known;  plus terms involving the first tangential derivatives of $\tilde c c^{-1}-1$ which we estimated  already. We can now deduce the desired estimate for $ {\partial (c^2\tilde c^{-2}-1)}/{\partial x^n}$ at $x^n=0$. 
  Differentiate again  \r{a3a} and \r{a3''} with respect to $x^n$ at $x^n=0$, to prove the estimate for $k=2$, etc. 
\end{proof}

\subsection{Local Interior Stability, proof of  Theorem~\ref{thm_stability}}
Set
\be{delta}
\delta := \|d^2-\tilde d^2\|_{C \big(\overline{\Gamma'\times \Gamma'}\big)}.
\ee
We use below interpolation estimates in the $C^k$ spaces, see, e.g., \cite{Triebel}.  If $m<k$, see \r{est}, we have
\be{d1'}
\|c-\tilde c\|_{C^m (U)} \le C\eps_0^{\mu}, 
\ee
with some $0<\mu= \mu(m)\le1$. Also,  by Theorem~\ref{thm_st_bd}, for any $m$, 
\be{d1}
\|d^2-\tilde d^2\|_{C^m \big(\overline{\Gamma'\times \Gamma'}\big)} \le C\delta^{\mu}, \quad 
\ee
with another $0<\mu= \mu(m)\le1$, if $k\gg 1$, under the a priori estimate \r{est} of the theorem. 
We will use the smallest $\mu$ above, and in the proof below, we will not specify the values of $k$ and $\mu$ which are guaranteed to work; even though this can be done. In principle, increasing $k$ in the theorem (assuming a stronger  a priori bound) increases $\mu$ (makes the H\"older stability estimate stronger), and vice-versa. 

We extend $c$ and $\tilde c$ outside $M$ by preserving the $\delta^\mu$-closeness on $\Gamma$. Extend $c$ in a smooth way first. 
As above, let $(x',x^n)$ be  boundary normal coordinates so that $x^n$ is the distance) in the metric $g_0$ to $\bo$, positive outside $M$. Given an integer $k$, let $E_k(\tilde c-c)$ be the truncated Taylor series of $\tilde c-c$ w.r.t.\  $x^n$ at $x^n=0$ of degree $m$, for $x^n>0$.

Let $0<-\rho_0\ll1$, $\mathsf{c}>0$ be such that the estimate in Proposition~\ref{pr_3.2} holds
in $\Omega_{\mathsf{c},\rho_0}$, with the choice of a boundary defining function $\rho=-x^n$. 
We can also assume that Proposition~\ref{pr_3.2}  holds for all sound speeds $c$ as in \r{est}, with $A$ fixed, and $\eps_0\ll1$. 
Set 
\[
\tilde c_1 = \left\{
\begin{array}{ll}
\tilde c, &x^n<0,\\
  c+ \chi(x^n/|\rho_0|)E_k(\tilde c-c), & x^n>0. 
\end{array}
\right.
\]
where $\chi(t)=1$ for $t<1/4$, $\chi(t)=0$ for $t\ge1/2$.  Then $\tilde c_1=\tilde c$ in $M$, and $\tilde c_1(x)=c(x)$  when $d(x,M)>|\rho_0|/2$. Moreover, for any small enough neighborhood $U$ of $p$, 
\be{d1a}
\| c-\tilde c_1\|_{C^m(U\setminus M)}\le C\|c-\tilde c\|_{C^m(\bo)}\le C\delta^\mu. 
\ee
We drop the subscript $1$  and denote the $\tilde c_1$ by $\tilde c$ below. 

Next, we compare $L$ and $\tilde L$ pushed to $\partial S B(x_0,r)$, see the Figure~\ref{fig:local_lens_rigidity_pic1}. Recall that in Section~\ref{sec_2.2}, we redefined $L$ and $\tilde L$ to act from $U_-$ to $U_+$, where $U_\pm\subset \partial_\pm  SB(x_0,r)$.

It is easy to see that fo $x$, $y$ near $x_0$, $d(x,y)/C\le \tilde d(x,y)\le Cd(x,y)$ with $C$ depending on $A$ in \r{est} only. 
Then writing $d-\tilde d = (d^2-\tilde d^2)/(d+\tilde d)$, we get from \r{d1}, 
\be{d2}
|\partial_{x',y'}^\alpha( d(x,y)-\tilde d(x,y))|\le C\delta^\mu /d(x,y)^{|\alpha|+1} \quad \text{on $\Gamma'\times\Gamma'$}.
\ee
We are only going to need this for $|\alpha|\le1$. 
Set
\[
\xi(x,y) = -\d_x d(x,y), \quad \eta(x,y) = \d_y d(x,y)
\]
which are just the tangent unit co-directions at $x$ and $y$ of the geodesic connecting  those two points.  
We define $\tilde\xi(x,y)$ and $\tilde \eta(x,y)$ in a similar way.  By the strict convexity,
\[
|\xi_n|\ge d(x,y)/C.
\]
Since $\xi$ and $\tilde \xi$ are unit, 
\[
\xi_n(x,y) = -\sqrt{c^{-2}(x)-g_0^{\alpha\beta}(x)\xi_\alpha\xi_\beta}, \quad \tilde \xi_n(x,y) = -\sqrt{\tilde c^{-2}(x)-g_0^{\alpha\beta}(x)\tilde \xi_\alpha \tilde \xi_\beta}. 
\]
Then, using \r{d2} with $|\alpha|=1$, we get
\[
|\xi_n(x,y) - \tilde \xi_n(x,y)| 
\le C\delta^\mu/ d(x,y)^{2}. 
\]
This, combined with \r{d2} for $\alpha=0$ yields
\be{d8}
| \xi(x,y) - \tilde \xi(x,y) |   \le C\delta^\mu/  d(x,y)^{2}. 
\ee
In the same way, we get
\be{d9}
\|\eta(x,y) - \tilde \eta(x,y)\|_{C^{k-1}}    \le C\delta^\mu/ d(x,y)^2. 
\ee

\begin{figure}[h!] 
  \centering
  \includegraphics[scale=0.75]{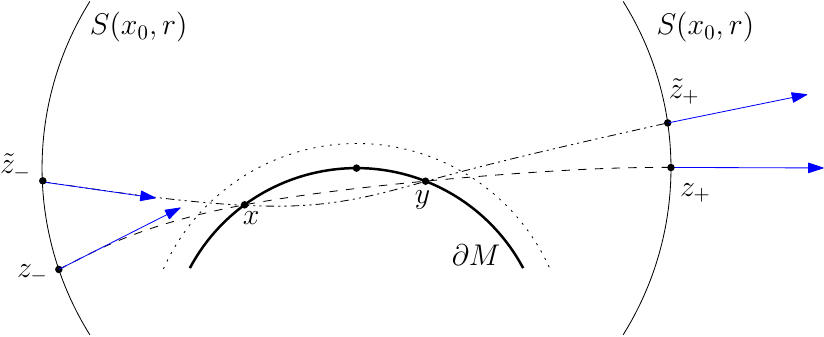}
\caption{Illustration to the proof of Theorem~\ref{thm_stability}. Here, $z_\pm$ and $\tilde z_\pm$ are actually points are codirections}
  \label{fig:local_lens_rigidity_pic3}
\end{figure}

For every $(x,y)\in\bo\times\bo$, both close to $x_0$, consider the shortest geodesics $\gamma_{[x,y]}$ and $\tilde \gamma_{[x,y]}$ connecting them, in the metrics $g$ and $\tilde g$, respectively. By the a priori estimate \r{est}, those two geodesics lie in a neighborhood of $x_0$ of size which shrinks to zero, as $d(x,x_0)$ and $d(y,x_0)$ tends to zero, uniformly for all $c$ and $\tilde c$ satisfying \r{est}. Next, we extend each one of them, lifted to the cotangent bundle (i.e., to the bicharacteristcs, called ''rays'' below), in both directions until they hit $\partial S B(x_0,r)$.  We denote by  $z_-= (x_-,  \xi_-)$   (respectively  $\tilde z_-= (\tilde x_-, \tilde \xi_-)$), the common point  with  $\partial_- S B(x_0,r)$, and by  $z_+=( y_+, \eta_+)$ (respectively  $\tilde z_+= (\tilde y_+,\tilde \eta_+)$), the  common point with  $\partial_+ S B(x_0,r)$, see Figure~\ref{fig:local_lens_rigidity_pic3}.

Let $\Sigma_-$ be the submanifold of $\partial_-SB(x_0,r)$ consisting of all $z$ there such that the bicharacteristic $Z(t,z)$ through it, in the metric $g$, is tangent to $\bo$ at some point close enough to $x_0$. Let $U_0$ be the connected neighborhood of it consisting of those point staying at distance (with respect to fixed coordinates) at most $O(\delta^{\mu'})$ with $0<\mu'<\mu$ which we chose below. Then the small but $\delta$-independent neighborhood $U_-$ of $\Sigma_-$, admits  the decomposition $U_- =  U_0+ U_{\rm in} + U_{\rm out}$, where $U_{\rm in}$ and $U_{\rm out}$ are disconnected components, such that the rays form  $U_{\rm in}$ in the metric $g$ hit $\bo$; and the   rays  from  $U_{\rm out}$ do not. All of them are extended until they hit  $\partial_+SB(x_0,r)$. We denote the corresponding travel times by $\ell(z)$, as in Section~\ref{sec_2.2} but now $\bo$ is replaced by the geodesic sphere $\partial B(x_0,r)$.

The rays issued from $U_{\rm out}$ would hit $\partial_+SB(x_0,r)$ and miss $\bo$, by definition. They will stay at distance at least $\delta^{\mu'}/C$ from $\bo$. The same is true for the rays  issued from $U_{\rm out}$  
related to the metric $\tilde g= \tilde c^{-2}g_0$ for $\delta\ll1$, because by \r{d1a}, the Hamiltonian vector fields of $g$ and $\tilde g$ are $O(\delta^\mu)$ close, and $\mu'<\mu$. In particular, the rays related to $\tilde g$ will not hit $\bo$ for $\delta\ll1$. Then
\be{d9'}
|Z(z,t)- \tilde Z(z,t)|\le C\delta^\mu, \quad z\in  U_{\rm out}
\ee
for $t$ over a compact interval. Recall that $Z(z,t)$ is the bicharacteristic issued from $z$.

The rays issued  from  $U_{\rm in}$ will hit $\bo$ at an angle at least $\delta^{\mu'/2}/C$, by the strict convexity. This is true, for $\delta\ll1$, for the rays related to $\tilde g$ as well by the closeness of the Hamiltonians (outside $M$) argument. 
The travel times $\tilde \ell(z)$ of the rays in the metric $\tilde g$ are  in general different than $\ell(z)$  
but the points  $\tilde Z(\ell(z),z)$ lie at distance to  $ Z(t,\ell(z))\in \partial_+SM$ which decreases with $\eps$ in \r{est}. Therefore, for $0<\eps\ll1$,   $\tilde Z(t,z)$ with  $z\in U_{\rm in}$  would enter $M$ and leave it for $t=\ell(z)$, and it will also leave the $O(\delta^\mu)$ neighborhood of  $M$ where we modified $\tilde c$. This is needed for the argument below. 

By \r{d8},  \r{d9},
\be{d10}
|z_--\tilde z_-|\le C\delta^\mu /d^2, \quad |z_+ -\tilde z_+|\le C\delta^\mu /d^2, z\in  U_{\rm in}.
\ee
Here, $d=d(x,y)$ with $(x,y)$ uniquely determined as the contacts of the ray trough $z$ with $\bo$. 
Also, by the closeness   of the Hamiltonians  outside $M$, and by \r{d9}, the lengths of the segments of $Z(z,t)$ and  $\tilde Z(\tilde z,t)$ from $z_-$ and $\tilde z_-$, respectively, to their common point $x$, differ by $O(\delta^\mu/d^2)$, see Figure~\ref{fig:local_lens_rigidity_pic3}. The segments in $M$ differ by $O(\delta^\mu/d)$, see \r{d2}. The remaining segments, after they exit $M$ and hit $S(x_0,r)$, differ by $O(\delta^\mu/d^2)$ by the above argument. Therefore,

\be{d11}
|\ell(z)-\tilde \ell(\tilde z)| \le C\delta^\mu/d^2,\quad  z\in  U_{\rm in}. 
\ee

By \r{1}, 
\be{d12}
\tilde Z(\ell(z_-),z_-) - Z(\ell(z_-),z_-)= 
\int_0^{\ell(z_-)} \frac{\partial\tilde Z}{\partial z} (t-s,Z (s,z_-))\big(\tilde V - V\big)(Z (s,z_-))\,\d s.
\ee
We replace next $\tilde Z(\ell(z_-),z_-)$ by $\tilde Z\big(\tilde \ell(\tilde z_-),\tilde z_-\big)=\tilde z_+$ modulo errors controlled by  \r{d10}, \r{d11},  so that we can get $z_+-z_-$ on the left; and use the fact that the latter satisfies \r{d10}. This would allow us to conclude that the integral is ``small'' in $U_{\rm in}$. 
 Indeed, since $Z(z,t)$ is a smooth function with uniformly bounded derivatives for $z$ and $t$ bounded, we get by \r{d10}, \r{d11},
\be{Zz}
\left|\tilde Z(\ell(z_-),z_-) - \tilde Z\big(\tilde \ell(\tilde  z_-),\tilde z_-\big)\right| 
\le  C\delta^\mu/d^2,\quad  z\in  U_{\rm in}. 
\ee
This and \r{d12} yield
\be{d13}
\begin{split}
\int_0^{\ell(z_-)}  \frac{\partial\tilde Z}{\partial z}& (t-s,Z (s,z_-))  \big(\tilde V - V\big)(Z (s,z_-))\,ds \\ &= \tilde z_+ - z_+ + O(\delta^\mu/d^2) = O(\delta^\mu/d^2),\quad  z\in  U_{\rm in},
\end{split}
\ee
where, in the last step, we used the second inequality in \r{d10}. The factor $d(x,y)$ can be replaced by the square root of the distance between $ U_{\rm in}$ and the tangent manifold $\Sigma_-$; and that root is $\delta^{\mu'/2}$ by construction. We choose $\mu'$ so that $\mu-2\mu'/2>\mu/2$ to obtain $O(\delta^{\mu/2}$) in \r{d13}. 

In $U_{\rm out}$,  the even better estimate $O(\delta^\mu)$ holds, without the $d^{-2}$ term, by \r{d9'}. In $U_0$, we used the argument explained before to extend the estimate from $U_{\rm out}$ to get an $O(\delta^\mu)+O(\delta^{\mu'}) = O(\delta^{\mu'}) $ estimate there. By interpolation, we get
\be{JJ}
\int_0^{\ell(z_-)} \frac{\partial\tilde Z}{\partial z} (t-s,Z (s,z_-))\big(\tilde V - V\big)(Z (s,z_-))\,\d s = O(\delta^{\mu/2}) \quad \text{in $C^1(U )$}. 
\ee 
Now,  by \r{1}, \r{2a},  and \r{3} we can write the left-hand side  of \r{JJ} in the form $\tilde J(f,\partial f)$, and complete the proof of Theorem~\ref{thm_st_bd} by applying Proposition~\ref{pr_3.2}.

\subsection{Global stability} 
For the purpose of the next theorem, we will extend and slightly generalize the foliation condition to compact submanifolds  $M_0$ of $M$. Let, as before, $\tilde M$ be a neighborhood of $M$, and extend $c$ smoothly there.  Note that the tilde over $M$ is not an indication that it is related to $\tilde c$. 

\begin{definition}\label{def_5.1}
Let $M_0\subset M$ be compact. We say that $M_0$ can be foliated by strictly convex hypersurfaces if there exists 
 a smooth function $\rho: \tilde M\to[0,\infty)$ which level sets $\Sigma_t=\rho^{-1}(t)$, $t\le T$
 with some $T>0$, restricted to $M_0$, are strictly convex viewed from $\rho^{-1}((0,t))$ for $g$;   $d\rho$ is non-zero on these level sets,   $\Sigma_0\cap M=\emptyset$, and $M_0\subset \rho^{-1}([0,T])$. 
\end{definition}

Note that this definition is not equivalent to Definition~\ref{def_1.1} when $M_0=M$ because in Definition~\ref{def_1.1}, we allow $M\setminus \rho^{-1}([0,T])$ to be non-empty (but with empty interior). Indeed, for uniqueness, proving $c=\tilde c$ outside such a set suffices since $c$ and $\tilde c$ are at least continuous. For stability however, it is convenient to assume that this set is empty.

We show next that the foliation condition implies non-trapping. 

\begin{proposition}\label{pr_nont}
If $M_0$ can be foliated by strictly convex hypersurfaces, then any maximal geodesic in $M_0$ is of finite length. 
\end{proposition}

\begin{proof}
Assume that $s\mapsto \gamma(s)$ exists for  $s\in (0,\infty)$ and stays in $M_0$;  and  in particular, it never reaches $\bo_0$ for $s>0$. The function  $f(s):= \rho\circ\gamma(s)$ cannot have more than one critical points as a consequence of the implication  $f'(s)=0 \Longrightarrow f''(s)<0$. By shifting the $s$ variable, we can always assume that the possible critical point is negative. 
Then $f$ is a strictly decreasing function for $s>0$.  On the other hand, it is bounded by below, so it has a limit $\hat\tau\ge0$, as $s\to\-\infty$, which is also its  infimum. 
By compactness, there exists a sequence $s_j\to \infty$ (we can start with $s_j=j$ and take a subsequence) so that $(x_j,v_j) := (\gamma(s), \dot \gamma(s))$ converges to some  $(\hat x,\hat v)\in SM_0$. Next, $ \rho(x_j) \searrow \hat\tau$. The limit $\hat v$ must be tangent to $\Sigma_{\hat\tau}$ at $\hat x$, because we can easily obtain a contradiction with the strict convexity if it is not.  
Now, by the strict convexity of $\Sigma_{\hat\tau}$ again, there exists $\delta>0$ so that $\gamma_{\hat x,\hat v}(s)$ would hit $\Sigma_{\hat \tau-\delta}$ for some positive time. This property is preserved under a small perturbation of  $(\hat x,\hat v)$; and therefore applies to $(x_j,v_j)$ for $j\gg1$.  This contradicts the choice of $\hat \tau$ however. 
\end{proof}

In particular, we get the following.

\begin{corollary}\label{cor_5.1}
Assume that $M=M_0\cup M_1$, where $M_0$ can be foliated by strictly convex hypersurfaces, and $M_1$ is non-trapping. Then $M$ is non-trapping as well. 
\end{corollary}

Note that $M_1$ can be a point, or a small neighborhood of a point, which happens if the level surfaces of $\rho=c>0$ are diffeomorphic to spheres but $\rho=0$ degenerates into a point. Another example is when $M_1$ is simple, see also the remarks following Theorem~\ref{thm_2}.

The $C^k$ norm below is defined in a fixed finite atlas of local coordinate charts.  In the same way we define $\dist(L,\tilde L)$ and its $C(D)$ norm: in any coordinate system we can just take the supremum of $L-\tilde L$ and then the maximum over all charts.  They can be defined in an invariant way, in principle but we do not do that for the sake of simplicity.

\begin{theorem}\label{thm_stab_global}
Assume that $M_0\subset M$ can  be foliated by strictly convex hypersurfaces for $g=c^{-2}g_0$. Let $D\subset \partial_- SM$ be  a neighborhood of 
the compact set of all $\beta\in \partial_- SM\cap \partial_- SM_0$ consisting of the initial points of all geodesics $\gamma_\beta$ tangent to the intersections  of the strictly convex hypersurfaces with $M_0$. 
 Then with $k$, $\mu$, $c_0$, $c$, $\tilde c$, $\eps_0$ and $A$ as in Theorem~\ref{thm_stability}, we have the stability estimate
\be{stab_global}
\|c-\tilde c\|_{C^2(M_0)}\le C\|\dist(L,\tilde L)\|_{C(D)}^\mu
\ee
for $c$, $\tilde c$ satisfying \r{est}. 
\end{theorem}

\begin{proof}
We will show first that the estimate in Theorem~\ref{thm_stability}, near a boundary point $p$, can be written in the form \r{stab_global}.
For $x$ and $y$ on $\bo$, both close to the fixed point $p$, let $\gamma_0 :[0,1]\mapsto \bo$ be the boundary geodesic connecting $x$ (for $t=0$) and $y$ (for $t=1$). Set $v(s) = \exp_x^{-1}\gamma_0(s)$, where $\exp$ is the interior exponential map. 
Then $(y,w)= L(x,v(1)/|v(1)|_g)$; set $(\tilde y,\tilde w)=\tilde  L(x,v(1)/|v(1)|_{\tilde g})  $. Set 
\be{Ld}
\delta=\|\dist(L,\tilde L)\|_{C(D)}
\ee
Then $|y-\tilde y|\le \delta$.  On the other hand, in the notation following \r{L},  given by $L(x,v)=(y,w)$, 
\be{w1}
\frac{\d}{\d s} d(x,\gamma_0(s))= \langle w'(s),\dot\gamma_0(s)\rangle_g,
\ee
 where $w'(s)$ is the tangential component of $w(s)$, the second component of $L(x, v(s)/|v(s)|_g)$. 
We also have 
\be{w2}
\frac{ \d}{\d s} \tilde d(x,\gamma_0(s))= \langle\tilde  w'(s),\dot\gamma_0(s)\rangle_{\tilde g},
\ee
 where $\tilde w(s)$ is the second component of $\tilde L(x, v(s)/|v(s)|_{\tilde g})$. Integrate the difference of \r{w1} and \r{w2} w.r.t.\ $s$ from $0$ to $1$,  and use \r{Ld} and the fact that $c=\tilde c$   on $\bo$ near $x$ to get $d(x,y)-\tilde d(x,\tilde y) = O(\delta)$ (we get $O(\delta)d(x,y)$, actually). Note that this also proves that $\ell-\tilde\ell = O(\delta)$ near $S_{x_0}M$. 
 We remark that in the inductive step below, we will only have that  $c$ and $\tilde c$ are $O(\delta^\mu)$ close on $\bo$ near $x$, instead of being equal. Then we would  get  $O(\delta^\mu)$ instead of $O(\delta)$.  
So we get now $d^2(x,y)-\tilde d^2(x,\tilde y) = O(\delta)$ as well. Since $d^2$ and $\tilde d^2$ are smooth and  $|y-\tilde y|\le \delta$, we get 
\[
d^2(x,y)-\tilde d^2(x,  y) = O(\delta).
\]
The constant in this estimate depends on the a priori bound $A$ in \r{est}. We can therefore apply Theorem~\ref{thm_st_bd} to estimate the jet of $c-\tilde c$ at $\bo$ of any finite order.  We then extend $c$ and $\tilde c$ in a neighborhood of $\bo_0\cap\bo$ so that  $c$ and $\tilde c$ are $O(\delta^\mu)$ there, as in \r{d1a}.

For any geodesic in $M_0$ issued from $D$, we extend it to $\tilde M$ and push the scattering relations $L$ and $\tilde L$ to $\partial \tilde M$. We show below that we have $\dist(L,\tilde L)= O(\delta^\mu)$ for the so extended lens relations, and a similar estimate for $\ell-\tilde\ell$. 
 We choose $U_0$,  $U_{\rm in}$ and $U_{\rm out}$ as the proof of Theorem~\ref{thm_st_bd} but this time $U_0$ is an $O(\delta^\mu)$ neighborhood of $\Sigma$ (rather than $O(\delta^{\mu'})$, because the singualar factor $1/d^2$ is not present anymore). 
For any $(x,v)\in U_{\rm in}$,  let $z_-$, $\tilde z_-$ be the points (in the phase space) of contact with the geodesic issued from $(x,v)$ for negative ``times'' with $\partial\tilde M$ w.r.t.\ the metrics $g$ and $\tilde g$, respectively. We define  $z_+$, $\tilde z_+$ in a similarly way, see Figure~\ref{fig:local_lens_rigidity_pic4}, where, in this case, $\Sigma_i=\bo$, and $\Omega_i= M_0$. 
Then $z_--\tilde z_-=O(\delta^\mu)$ (the difference makes sense in local coordinates only) because we have the same for $c-\tilde c$ in $C^2$. We have an $O(\delta)$ bound of $y-\tilde y$, and of the difference of the corresponding directions there, by the closeness of the scattering relations. Next, $z_+-\tilde z_+=O(\delta^\mu)$ as well. Therefore, $\dist(L(z_-),\tilde L(\tilde z_-))= O(\delta^\mu)$. 
Since $L$ and $\tilde L$ are smooth with uniformly bounded first derivatives, we also get $\dist(L(z_-),\tilde L(z_-))= O(\delta^\mu)$, for all $z_-\in U_{\rm in}$. We also have $\ell(z_-)-\tilde \ell (z_-) = O(\delta^\mu)$ as shown in the paragraph following \r{w2}.  
In this step, we are no longer working near a fixed point on $\bo$, and for almost tangential directions, the requirement for $U_{\rm in}$ is that the geodesic issued from it are in the set required by the theorem but that set is a distance at least $\delta^\mu$ from the tangential set $\Sigma$. In $U_{\rm out}$ and $U_{0}$, we argue as above. 

Therefore, we have now that $\|\dist(L,\tilde L)\|_{C(D_0)}= O(\delta^\mu)$, where $L$,  $\tilde L$ and $D_0$ are defined as $L$,  $\tilde L$ and $D$ but with $\bo$ replaced by $\partial \tilde M$. We have the same for $\ell-\tilde \ell$ as well.  The advantage of this is that the geodesics issued from $D_0$ hitting $\supp(c-\tilde c) $ are never tangent to $\partial \tilde M$, and we actually have a uniform lower bound on the angle  they make with $\partial \tilde M$. Thus we reduced the analysis to the case when $c-\tilde c$ is a priori known to be supported in the interior of $M$.

\begin{figure}[h!] 
  \centering
  \includegraphics[scale=0.75]{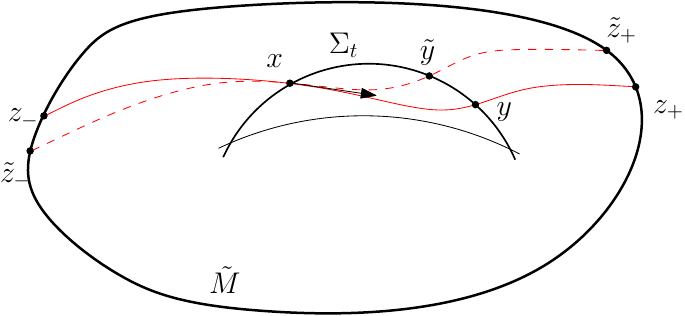}
\caption{Illustration to the proof of Theorem~\ref{thm_stab_global}.  
}
  \label{fig:local_lens_rigidity_pic4}
\end{figure}

We then apply a layer stripping argument finite number of times, see also the introduction in \cite{UV:Local}. For each $\Sigma_t$, assuming the estimate outside $\Sigma_t$, we can choose an appropriate lens shaped domain near each point on $\Sigma_t$, and artificial boundary close enough to it,  see  Figure~\ref{fig:local_lens_rigidity_pic4}. The $O(\delta^\mu)$ closeness of $L$ and $\tilde L$, and $\ell$ and $\tilde \ell$ implies the same for $c$ and $\tilde c$ by the arguments using \r{JJ}. We do this for a finite number of points on $\Sigma_t$ by a compactness argument to prove the estimate in $\rho^{-1}(0,t+\epsilon)$ for some
 $\epsilon>0$.  We then cover $M_0$ with a finite set  $\Sigma_i := \rho^{-1}(t_i)$, $M_i:=  \rho^{-1}(-\infty, t_i)$ and finish the proof with finitely many steps. 
\end{proof}

\section{Herglotz and Wieckert \& Zoeppritz speeds are lens rigid}\label{section_6}

We revisit the 
Herglotz \cite{Herglotz} and Wieckert and Zoeppritz \cite{WZ} class of speeds. Let $M = B(0,R)$, $R>0$ be the ball in $\R^n$, $n\ge3$ centered at the origin with radius $R>0$. The background metric $g_0$ in this section is the Euclidean one. 
Let $0<c(x)$ be smooth in $B(0,R)$. 
Assume that $c$ satisfies the Herglotz and Wieckert \& Zoeppritz condition
\be{HWZ}
\frac{\partial}{\partial r}\frac{r}{c(r)}>0, \quad \text{for $0<r=|x|\le R$},
\ee
where $\frac{\partial}{\partial r} = \frac{x}{|x|}\cdot\partial_x$ is the radial derivative. 
We do not assume that $c$ is radial, i.e., that it depends on $r=|x|$ only. We show below that \r{HWZ} is in fact a foliation condition. 

\begin{proposition}\label{pr_HWZ}
The Herglotz and Wieckert \& Zoeppritz condition \r{HWZ} is equivalent to the the condition that the Euclidean spheres $S_r=\{|x|=r\}$ are strictly convex in the metric $c^{-2} \d x^2$ for $0<r\le R$. 
\end{proposition}

\begin{proof} This proposition is essentially proved in \cite[Proposition~7.1]{S-U-InsideOut10}. 
We will show first that the strict convexity condition of an oriented hypersurface $S$ (positive second fundamental form) is equivalent to the following: (A) if $G$ is the generator of the geodesic flow, and if $\rho$ is a defining function of $S$ positive on the ``positive'' side of $S$, then $G^2 \rho>0$ when $G\rho=\rho =0$ on a non-zero the energy level. 
In semigeodesic coordinates $(x',x^n)$ with $x^n>0$  on the ``positive'' side of $S$, we have $\rho = x^n \rho'$ with $\rho'>0$ on $S$. Since
\[
G = \xi^j \partial_{x^j} - \Gamma_{ij}^k \xi^i\xi^j\partial_{\xi_k},
\]
we have $G\rho = \xi^n+ x^n G\rho'$. On $S$, $G\rho= \xi^n$; therefore, (A) can be reformulated as  follows: $x^n=\xi^n=0$ and $\xi\not=0$ imply $G^2\rho>0$. Differentiating $G\rho$ again, we see that under those conditions, $G^2\rho = - \Gamma_{\alpha\beta}^n\xi^\alpha \xi^\beta$, where the Greek indices run from $1$ to $n-1$. This is well known and follows directly form the definition to be the second fundamental form corresponding to the orientation $x^n>0$. Then strict convexity, viewed from $\rho>0$ is equivalent to a negative second fundamental form. For manifolds with boundary we study,  $\rho>0$ in the interior of $M$, therefore convexity viewed from the exterior means  a positive second fundamental form.

We can take $\rho=r^2/2-r_0^2/2=|x|^2/2-r_0^2/2$ as a defining function of the sphere $S_{r_0}$ for $x\not=0$. Then   $G\rho=\xi\cdot x$. Next, $G^2\rho = |\xi|^2- \Gamma_{ij}^k\xi^i\xi^j x^k$ and on the energy level $c^{-2}|\xi|^2=1$, we have  $G^2\rho = c^2- \Gamma_{ij}^k\xi^i\xi^j x^k$. Here and below, we still have summation w.r.t.\ $k$ even though both indices are upper. A direct computation shows that 
\[
\Gamma_{ij}^k = \frac12 c^2\left(\delta_{j}^k\partial_{x^i} + \delta_{i}^k\partial_{x^j} -  \delta_{ij}\partial_{x^k}\right)c^{-2},
\]
with $c=c(|x|)$, therefore,
\[
\Gamma_{ij}^k\xi^i\xi^j x^k=- c^{-1}\left(2\xi^kx^k\xi^i \partial_{x^i}-|\xi|^2 x^k \partial_{x^k}\right) c = -c^{-1}\left( 2(\xi\cdot x) \xi\cdot\partial_x - |\xi|^2 x\cdot\partial_x  \right)c.
\]
Therefore, on the unit energy level and for $G\rho=0$, strict convexity is equivalent to 
\[
c -x\cdot\partial_x c>0,
\]
i.e., $r\partial_r c<c$, which is equivalent to \r{HWZ}. Note that the computations in \cite[Proposition~7.1]{S-U-InsideOut10} are done in the cotangent bundle and are somewhat shorter. 
\end{proof}

\begin{corollary}\label{cor_HWZ}
Speeds satisfying \r{HWZ} are lens rigid in $B(0,R)$ in the sense of Theorem~\ref{thm_2}. 
\end{corollary}

Note that we only require $c$ in Theorem~\ref{thm_2} to satisfy \r{HWZ}. Then any other speed  $\tilde c$ for which $\partial B(0,R)$ is still strictly convex and $\tilde c=c$ on $\partial B(0,R)$ with the same scattering relation is equal to $c$. This extends the Herglotz \cite{Herglotz} and the Wieckert \& Zoeppritz \cite{WZ} results to not necessarily  radial speeds $c$ satisfying \r{HWZ}. 

If \r{HWZ} holds in the shell $R_0<|x|<R$ only, with some $0<R_0<R$, then we get lens rigidity in the shell and we only need to use the scattering relation for geodesics staying in it. The speed $c$ does not even need to be defined in $|x|<R_0$. We will skip the details.



\end{document}